\documentclass[11pt]{article}

\usepackage{epsfig}
\usepackage{latexsym}
\usepackage{url}
\usepackage{float}
\usepackage[hypertexnames=false,hyperfootnotes=false]{hyperref}
\usepackage{texnansi}
\usepackage{color}
\usepackage{tikz}
\usepackage{subfigure}
\usepackage{afterpage}
\usepackage{enumerate}
\usepackage{algorithm}
\usepackage{algorithmic}
\usepackage[normalem]{ulem}
\usepackage{lmodern}
\usepackage{sectsty}
\usepackage{ifthen}
\usepackage[leqno]{amsmath}
\usepackage{amsthm}
\usepackage{amssymb}
\usepackage{amsfonts}
\usepackage[margin=10pt,font=small,labelfont=bf]{caption}
\usepackage{natbib}
\newcommand{\field}[1]{\ensuremath{\mathbb{#1}}}
\newcommand{\R}{\ensuremath{\field{R}}} 
\newcommand{\tends}{\ensuremath{\rightarrow}} 
\newcommand{\E}{\ensuremath{\mathsf{E}}} 
\newcommand{\defeq}{\ensuremath{\triangleq}}
\renewcommand{\Re}{\ensuremath{\R}} 

\DeclareMathOperator{\Var}{Var}

\DeclareMathOperator*{\argmax}{\text{argmax}}

\renewcommand\max{\ensuremath{\mathop{\text{max}}\limits}}

\renewcommand\lim{\ensuremath{\mathop{\text{lim}}\limits}}

\renewcommand\exp{\ensuremath{\mathop{\text{exp}}\limits}}

\newtheoremstyle{thm-sf}{}{}{\itshape}{}{\sffamily\bfseries}{.}{ }{}
\theoremstyle{thm-sf}

\newtheorem{definition}{Definition}

\newtheorem{theorem}{Theorem}

\newtheorem{lemma}{Lemma}

\usetikzlibrary{arrows,patterns,plotmarks}
\tikzstyle{every picture} += [>=stealth]
\allsectionsfont{\sffamily}
\makeatletter 
\def\@seccntformat#1{\csname the#1\endcsname.\quad} 
\makeatother

\floatstyle{ruled}
\newcommand{\emailhref}[1]{\href{mailto:#1}{\tt #1}} 
\provideboolean{fastcompile}
\newcommand{\hidefastcompile}[1]{\ifthenelse{\boolean{fastcompile}}{}{#1}}

\usepackage{setspace}
\newcommand{\algeqstart}{\vspace*{-0.5\baselineskip}}
\newcommand{\algeqend}{\vspace*{-1.25\baselineskip}}

\setcitestyle{authoryear,round,semicolon,aysep={,},yysep={,},notesep={, }}
\newtheorem*{th:Gaussian}{Theorem~\ref{th:Gaussian}}
\newtheorem*{th:tqd-aqd}{Theorem~\ref{th:tqd-aqd}}
\newtheorem*{th:sufficientcondition}{Theorem~\ref{th:sufficientcondition}}
\newtheorem*{th:linear-pi-psi}{Theorem~\ref{th:linear-pi-psi}}

\title{\sffamily\bfseries
  Strategic Execution \\
  in the Presence of an Uninformed Arbitrageur
}
\author{
{\sffamily Ciamac C. Moallemi} \\
Graduate School of Business \\
Columbia University \\
email: \emailhref{ciamac@gsb.columbia.edu} \\
\and
{\sffamily Beomsoo Park} \\
Electrical Engineering \\
Stanford University \\
email: \emailhref{beomsoo@stanford.edu} \\
\and
{\sffamily Benjamin Van Roy} \\
Management Science \& Engineering \\
Electrical Engineering \\
Stanford University \\ 
email: \emailhref{bvr@stanford.edu}\\
}
\date{January 18, 2008 \\
  Revised: February 28, 2009}

\begin{document}

\maketitle

\singlespacing

\begin{abstract}
  We consider a trader who aims to liquidate a large position in the
  presence of an arbitrageur who hopes to profit from the trader's
  activity.  The arbitrageur is uncertain about the trader's position
  and learns from observed price fluctuations.  This is a dynamic game
  with asymmetric information.  We present an algorithm for computing
  perfect Bayesian equilibrium behavior and conduct numerical
  experiments.  Our results demonstrate that the trader's strategy
  differs significantly from one that would be optimal in the absence
  of the arbitrageur. In particular, the trader must balance the
  conflicting desires of minimizing price impact and minimizing
  information that is signaled through trading. Accounting for
  information signaling and the presence of strategic adversaries can
  greatly reduce execution costs.
 \end{abstract}

\doublespacing

\section{Introduction}

When buying or selling securities, value is lost through execution
costs such as exchange fees, commissions, bid-ask spreads, and price
impact.  The latter can be dramatic and typically dominates other
sources of execution cost when trading large blocks, when the security
is thinly traded, or when there is an urgent demand for liquidity.
Execution algorithms aim to reduce price impact by partitioning the
quantity to be traded and placing trades sequentially.  Growing
recognition for the importance of execution has fueled an academic
literature on the topic as well as the formation of specialized groups
at investment banks and other organizations to offer execution
services.

Optimal execution algorithms have been developed for a number of
models.  In the base model of \citet{BertsimasLo98}, a
stock price nominally follows a discrete-time random walk and the
market impact of a trade is permanent and linear in trade size.  The
authors establish that expected cost is minimized by an
equipartitioning policy. This policy trades equal amounts over 
time increments within the trading horizon.  Further developments have
led to optimal execution algorithms for models that incorporate price
predictions \citep{BertsimasLo98}, bid-ask spreads and resilience
\citep{ObizhaevaWang05,AlfonsiSchiedSchulz07a}, nonlinear price impact
models \citep{Almgren03,AlfonsiSchiedSchulz07b}, and risk aversion
\citep{SubramanianJarrow01,AlmgrenChriss00,Dubil02,HubermanStanzl05,Engle06,Hora06,Almgren06,SchiedSchonenborn07,Lorenz08}.

The aforementioned results offer insight into how one should partition
a block and sequence trades under various assumptions about market
dynamics and objectives.  The resulting algorithms, however, are
unrealistic in that they exhibit predictable behavior. Such
predictable behavior allows strategic adversaries, which we call
arbitrageurs, to ``front-run'' trades and profit at the expense of
increased execution cost.  For example, consider liquidating a large
block by an equipartitioning policy which sells an equal amount during
each minute of a trading day.  Trades early in the day generate
abnormal price movements, allowing an observing arbitrageur to
anticipate further liquidation.  If the arbitrageur sells short and
closes his position at the end of the day, he profits from expected
price decreases.  The arbitrageur's actions amplify price impact and
therefore increase execution costs.

Several recent papers study game-theoretic models of execution in the
presence of strategic arbitrageurs
\citep{BrunnermeierPedersen05,CarlinLoboViswanathan07,SchonenbornSchied07}.
However, these models involve games with symmetric information, in
which arbitrageurs know the position to be liquidated. In more
realistic scenarios, this information would be the private knowledge
of the trader, and the arbitrageurs would make inferences as to the
trader's position based on observed market activity.


This type of information asymmetry is central to effective
execution. The fact that his position is unknown to others allows the
trader to greatly reduce execution costs. But to do so requires
the deliberate management of ``information leakage'', or the signals that
are transmitted via trading activity. Further, the desire to minimize
information signaling may be at odds with the desire to minimize price
impact.  A model through which such signaling can be studied must
account for uncertainty among arbitrageurs and their ability to learn
from observed price fluctuations.  In this paper we formulate and
study a simple model which we believe to be the first that meets this
requirement.

The contributions of this paper are as follows:
\begin{enumerate}
  \setlength{\itemsep}{0pt}
\item We formulate the optimal execution problem as a dynamic game
  with asymmetric information. This game involves a trader and a
  single arbitrageur. Both agents are risk neutral, and market
  dynamics evolve according to a linear permanent price impact
  model. The trader seeks to liquidate his position in a finite time
  horizon. The arbitrageur attempts to infer the position of the
  trader by observing market price movements, and seeks to exploit
  this information for profit.

\item We develop an algorithm that computes perfect Bayesian
  equilibrium behavior. 

\item We demonstrate that the associated equilibrium strategies take
  on a simple structure: Trades placed by the trader are linear in the
  trader's position, the arbitrageur's position and the arbitrageur's
  expectation of the trader's position. Trades placed by the
  arbitrageur are linear in the arbitrageur's position and his
  expectation of the trader's position. Equilibrium policies
  depend on the time horizon and a parameter that we
  call the ``relative volume''. This parameter captures the magnitude
  of the per-period activity of the trader relative to the exogenous
  fluctuations of the market.

\item We present computational results that make several points about
  perfect Bayesian equilibrium in our model: 
  \begin{enumerate}
    \setlength{\itemsep}{0pt}
  \item In the presence of adversaries, there are significant
    potential benefits to employing perfect Bayesian equilibrium
    strategies. 
  \item Unlike strategies proposed based on prior models in the
    literature, which exhibit deterministic sequences of trades,
    trades in a perfect Bayesian equilibrium adaptively respond to price
    fluctuations; the trader leverages these random outcomes to conceal
    his activity.
  \item When the relative volume of the trader's activity is low, in
    equilibrium, the trader can ignore the presence of the arbitrageur
    and will equipartition to minimize price impact. Alternatively,
    when the relative volume is high, the trader will concentrate his
    trading activity in a short time interval so as to minimize
    signaling.
  \item The presence of the arbitrageur leads to a spill-over effect.
    That is, the trader's expected loss due to the arbitrageur's
    presence is larger than the expected profit of the
    arbitrageur. Hence, other market participants benefit from the
    arbitrageur's activity.
  \end{enumerate}

\item We discuss how the basic model presented can be can be extended
  to incorporate a number of additional features, such as transient
  price impact and risk aversion.
\end{enumerate}


Solving for perfect Bayesian equilibrium in dynamic games with asymmetric
information is notoriously difficult.  What facilitates effective computation
in our model is that, in equilibrium, each agent solves a tractable
linear-quadratic Gaussian control problem.  Similar approaches based on
linear-quadratic Gaussian control have previously been used to analyze
equilibrium behavior of traders with private information. This line of work
begins with the seminal paper of \cite{Kyle85}, and includes many subsequent
papers \citep[e.g.,][]{Foster94,Foster96,Vayanos01}.  Among these
contributions, \citet{Foster94} come closest to the model and method we
propose.  In the model of that paper, there are two strategic traders, many
``noise'' traders, and a market maker.  The strategic traders possess
information that is not initially reflected in market prices.  One trader
knows more than the other.  The more informed trader adapts trades to maximize
his expected payoff, and this entails controlling how his private information
is revealed through price fluctuations.  This model parallels ours if we think
of the arbitrageur as the less informed trader.  However, in our model there
is no private information about future dividends but instead uncertainty about
the size of the position to be liquidated.  Further, in the model of
\citet{Foster94}, trades influence prices because the market maker tries to
infer the traders' private information whereas, in our setting, there is an
exogenously specified price impact model.  The algorithm we develop bears some
similarity to that of \citet{Foster94}, but requires new features designed to
address differences in our model.


The remainder of this paper is organized as follows.  The next section
presents our problem formulation.  Section~\ref{se:dp-analysis}
discusses how perfect Bayesian equilibrium in this model is
characterized by a dynamic program.  A practical algorithm for
computing perfect Bayesian equilibrium behavior is developed in
Section~\ref{se:algorithm}.  This algorithm is applied in
computational studies, for which results are presented and interpreted
in Section~\ref{se:computation}.  Several extensions of this model are
discussed in Section~\ref{se:extensions}.  Finally,
Section~\ref{se:conclusion} makes some closing remarks and suggests
directions for future work.  Proofs of all theoretical results are
presented in the appendices.

\section{Problem Formulation}\label{se:model}

In this section, the optimal execution problem is formulated as a game
of asymmetric information. Our formulation makes a number of
simplifying assumptions and we omit several factors that are important
in the practical implementation of execution strategies, for example,
transient price impact and risk aversion. Our goal here is to
highlight the strategic and informational aspects of execution in a
streamlined fashion. However, these assumptions are discussed in more
detail and a number of extensions of this basic model are presented in
Section~\ref{se:extensions}.

\subsection{Game Structure}

Consider a game that evolves over a finite horizon in discrete time
steps $t=0,\ldots,T+1$.  There are two players: a trader and an
arbitrageur.  The trader begins with a position $x_0 \in \Re$ in a
stock, which he must liquidate by time $T$. Denote his position at
each time $t$ by $x_t$, and thus require that $x_t=0$ for $t \geq
T$. The arbitrageur begins with a position $y_0$.  Denote his
position at each time $t$ by $y_t$.  In general, the arbitrageur has
additional flexibility and will not be limited to the same time
horizon as the trader. For simplicity, this flexibility is modelled by
assuming that the arbitrageur has one additional period of trading
activity. In other words, though we do require that $y_{T+1}=0$,
we do not require that $y_T=0$. This assumption
will be revisited in Section~\ref{se:trading-horizon}.

\subsection{Price Dynamics}

Denote the price of the stock at time $t$ by $p_t$. This price evolves
according to the permanent linear price impact model given by
\begin{equation}\label{eq:p-evolve}
p_{t} 
= p_{t-1} + \Delta p_{t}
= p_{t-1} + \lambda (u_t + v_t) + \epsilon_t.
\end{equation}
Here, $\lambda > 0$ is a parameter that reflects the sensitivity of
prices to trade size, and $u_t$ and $v_t$ are, respectively, the
quantities of stock purchased by the trader and the arbitrageur at
time $t$. Note that, given the horizon of the trader, $u_{T+1} \defeq
0$.  The positions evolve according to
\[
x_t = x_{t-1} + u_t, 
\quad \text{and} \quad y_t = y_{t-1} + v_t.
\]

The sequence $\{\epsilon_t\}$ is a normally distributed IID process with
$\epsilon_t \sim N(0,\sigma_\epsilon^2)$, for some $\sigma_\epsilon>0$. This
noise sequence represents the random and exogenous fluctuations of market
prices. We assume that the trading decisions $u_t$ and $v_t$ are made at time
$t-1$, and executed at the price $p_t$ at time $t$. Note that there is no
drift term in the price evolution equation \eqref{eq:p-evolve}. In the
intraday horizon of typical optimal execution problems, this is usually a
reasonable assumption. This assumption will be revisited in
Section~\ref{se:price-impact}. Further, the price impact in
\eqref{eq:p-evolve} is permanent in the sense that it is long-lived relative
to the length of the time horizon $T$. In Section~\ref{se:price-impact} we
will allow for transient price impact as well.imp

\subsection{Information Structure}

The information structure of the game is as follows. The dynamics of
the game (in particular, the parameters $\lambda$ and
$\sigma_\epsilon$) and the time horizon $T$ are mutually known. From
the perspective of the arbitrageur, the initial position $x_0$ of the
trader is unknown. Further, the trader's actions $u_t$ are not
directly observed. However, the arbitrageur begins with a prior
distribution $\phi_0$ on the trader's initial position $x_0$. As the
game evolves over time, the arbitrageur observes the price change
$\Delta p_t$ at each time $t$.  The arbitrageur updates his beliefs
based on these price movements, at any time $t$ maintaining a
posterior distribution $\phi_t$ of the trader's current position
$x_t$, based on his observation of the history of the game up to and
including time $t$.

From the trader's perspective, it is assumed that everything is known.
This is motivated by the fact that the arbitrageur's initial position
$y_0$ will typically be zero and the trader can go through the same
inference process as the arbitrageur to arrive at the prior
distribution $\phi_0$. Given a prescribed policy for the arbitrageur
(for example, in equilibrium), the trader can subsequently reconstruct
the arbitrageur's positions and beliefs over time, given the public
observations of market price movements. We do make the assumption,
however, that any deviations on the part of the arbitrageur from his
prescribed policy will not mislead the trader. In our context, this
assumption is important for tractability. We discuss the situation
where this assumption is relaxed, and the trader does not have perfect
knowledge of the arbitrageur's positions and beliefs, in
Section~\ref{se:conclusion}.

\subsection{Policies}

The trader's purchases are governed by a policy, which is a sequence
of functions $\pi = \{\pi_1,\ldots,\pi_T\}$. Each function $\pi_{t+1}$
maps $x_{t}$, $y_{t}$, and $\phi_{t}$, to a decision $u_{t+1}$ at time
$t$. Similarly, the arbitrageur follows a policy $\psi =
\{\psi_1,\ldots,\psi_{T+1}\}$.  Each function $\psi_{t+1}$ maps $y_t$
and $\phi_t$ to a decision $v_{t+1}$ made at time $t$. Since policies
for the trader and arbitrageur must result in liquidation, we require
that $\pi_{T}(x_{T-1},y_{T-1},\phi_{T-1}) = -x_{T-1}$ and
$\psi_{T+1}(y_{T},\phi_{T}) = -y_{T}$.  Denote the set of trader
policies by $\Pi$ and the set of arbitrageur policies by $\Psi$.

Note that implicit in the above description is the restriction to
policies that are Markovian in the following sense: the state of the
game at time $t$ is summarized for the trader and arbitrageur by the
tuples $(x_{t},y_{t},\phi_t)$ and $(y_t,\phi_t)$, respectively, and
each player's action is only a function of his state. Further,
the policies are pure strategies in the sense that, as a function of
the player's state, the actions are deterministic. In general, one may
wish to consider policies which determine actions as a function of the
entire history of the game up to a given time, and allow randomization
over the choice of action. Our assumptions will exclude equilibria
from this more general class. However, it will be the case that for
the equilibria that we do find, arbitrary deviations that are history
dependent and/or randomized will not be profitable.

If the arbitrageur applies an action $v_t$ and assumes the trader uses
a policy $\hat{\pi} \in \Pi$, then upon observation of $\Delta p_t$ at
time $t$, the arbitrageur's beliefs are updated in a Bayesian fashion
according to
\begin{equation}\label{eq:bayes-update}
\phi_t(S) = \Pr\big(x_t \in S \ |\ \phi_{t-1}, y_{t-1}, 
\lambda(\hat{\pi}_t(x_{t-1}, y_{t-1}, \phi_{t-1}) + v_t) 
+ \epsilon_t = \Delta p_t\big),
\end{equation}
for all measurable sets $S\subset\R$.  Note that $\Delta p_t$ here is
an observed numerical value which could have resulted from a trader
action $u_t \neq \hat{\pi}_t(x_{t-1}, y_{t-1}, \phi_{t-1})$.  As such,
the trader is capable of misleading the arbitrageur to distort his
posterior distribution $\phi_t$.

\subsection{Objectives}

Assume that both the trader and the arbitrageur are risk neutral and
seek to maximize their expected profits (this assumption will be
revisited in Section~\ref{se:risk-aversion}).  Profit is computed
according to the change of book value, which is the sum of a player's
cash position and asset position, valued at the prevailing market
price.  Hence, the profits generated by the trader and arbitrageur
between time $t$ and time $t+1$ are, respectively,
\[
p_{t+1} x_{t+1} - p_{t+1} u_{t+1} - p_t x_t = \Delta p_{t+1} x_t,
\quad\text{and}\quad
p_{t+1} y_{t+1} - p_{t+1} v_{t+1} - p_t y_t = \Delta p_{t+1} y_t.
\]

If the trader uses policy $\pi$ and the arbitrageur uses policy $\psi$
and assumes the trader uses policy $\hat{\pi}$, the trader expects
profits
\[
U_t^{\pi,(\psi,\hat{\pi})}(x_t,y_t,\phi_t) 
\defeq \E^{\pi,(\psi,\hat{\pi})}\left[\left.
\sum_{\tau=t}^{T-1} \Delta p_{\tau+1} x_\tau\ \right|\ x_t,y_t, \phi_t \right],
\]
over times $\tau = t+1, \ldots, T$.  Here, the superscripts indicate
that trades are executed based on $\pi$ and $\psi$, while beliefs are
updated based on $\hat{\pi}$.
Similarly, the arbitrageur expects profits
\[
V_t^{(\psi,\hat{\pi}),\pi}(y_t,\phi_t) 
\defeq  \E^{\pi,(\psi,\hat{\pi})}\left[\left.
\sum_{\tau=t}^{T} \Delta p_{\tau+1} y_\tau\ \right|\ y_t, \phi_t \right],
\]
over times $\tau = t+1, \ldots, T+1$. Here, the conditioning in the
expectation implicitly assumes that $x_t$ is distributed according to
$\phi_t$.

Note that $-U_t^{\pi,(\psi,\hat{\pi})}(x_0,y_0,\phi_0)$ is the
trader's expected execution cost.  For practical choices of $\pi$,
$\psi$, and $\hat{\pi}$, we expect this quantity to be positive since
the trader is likely to sell his shares for less than the initial
price.  To compress notation, for any $\pi$, $\psi$, and $t$, let
\[
U_t^{\pi,\psi}  \defeq U_t^{\pi,(\psi,\pi)},  
\quad \text{and} \quad 
V_t^{\psi,\pi}  \defeq V_t^{(\psi,\pi),\pi}.
\]

\subsection{Equilibrium Concept}

As a solution concept, we consider perfect Bayesian equilibrium
\citep{FudenbergTirole91}. This is a refinement of Nash equilibrium
that rules out implausible outcomes by requiring subgame perfection
and consistency with Bayesian belief updates. In particular, a policy $\pi \in \Pi$ is a {\em best response} to
$(\psi,\hat{\pi}) \in \Psi \times \Pi$ if
\begin{equation}
\label{eq:trader-response}
U_t^{\pi, (\psi,\hat{\pi})}(x_t,y_t,\phi_t) = \max_{\pi' \in \Pi}\ 
U_t^{\pi', (\psi,\hat{\pi})}(x_t,y_t,\phi_t),
\end{equation}
for all $t$, $x_t$, $y_t$, and $\phi_t$.  Similarly, a policy $\psi
\in \Psi$ is a {\em best response} to $\pi \in \Pi$ if
\begin{equation}
\label{eq:arbitrageur-response}
V_t^{\psi,\pi}(y_t,\phi_t) = \max_{\psi' \in \Psi}\ 
V_t^{\psi',\pi}(y_t,\phi_t),
\end{equation}
for all $t$, $y_t$, and $\phi_t$. We define perfect Bayesian
equilibrium, specialized to our context, as follows:
\begin{definition}
  A {\bf perfect Bayesian equilibrium} (PBE) is a pair of policies
  $(\pi^*,\psi^*) \in \Pi \times \Psi$ such that:
\begin{enumerate}
\item $\pi^*$ is a best response to $(\psi^*,\pi^*)$;
\item $\psi^*$ is a best response to $\pi^*$.
\end{enumerate}
\end{definition}

In a PBE, each player's action at time $t$ depends on positions $x_t$
and/or $y_t$ and the belief distribution $\phi_t$.  These arguments,
especially the distribution, make computation and representation of a
PBE challenging.  We will settle for a more modest goal.  We compute
policy actions only for cases where $\phi_t$ is Gaussian.  When the
initial distribution $\phi_0$ is Gaussian and players employ these PBE
policies, we require that subsequent belief distributions
$\phi_t$ determined by Bayes' rule \eqref{eq:bayes-update} also be
Gaussian.  As such, computation of PBE policies over the restricted
domain of Gaussian distributions is sufficient to characterize
equilibrium behavior given any initial conditions involving a Gaussian
prior.  To formalize our approach, we now define a solution concept.
\begin{definition}
  A policy $\pi \in \Pi$ (or $\psi \in \Psi$) is a {\bf Gaussian best
    response} to $(\psi,\hat{\pi}) \in \Psi \times \Pi$ (or $\pi \in
  \Pi$) if \eqref{eq:trader-response} (or
  \eqref{eq:arbitrageur-response}) holds for all $t$, $x_t$, $y_t$,
  and Gaussian $\phi_t$.  A {\bf Gaussian perfect Bayesian
    equilibrium} is a pair $(\pi^*,\psi^*) \in \Pi \times \Psi$ of
  policies such that
  \begin{enumerate}
  \item $\pi^*$ is a Gaussian best response to $(\psi^*,\pi^*)$;
  \item $\psi^*$ is a Gaussian best response to $\pi^*$;
  \item if $\phi_0$ is Gaussian and arbitrageur assumes the trader
    uses $\pi^*$ then, independent of the true actions of the trader,
    the beliefs $\phi_1, \ldots,\phi_{T-1}$ are Gaussian.
  \end{enumerate}
\end{definition}
Note that when Gaussian PBE policies are used and the prior $\phi_0$
is Gaussian, the system behavior is indistinguishable from that of a
PBE since the policies produce actions that concur with PBE policies
at all states that are visited.

Given a belief distribution $\phi_t$, define the quantities
\[
\mu_t  \defeq \E[ x_t\ |\ \phi_t], \qquad
\sigma_t^2 \defeq \E\left[\left.(x_t - \mu_t)^2\ \right| \phi_t\right], \quad\text{and}\quad
\rho_t \defeq \lambda \sigma_t /
\sigma_\epsilon.
\]
Since $\lambda$ and $\sigma_\epsilon$ are constants, $\rho_t$ is
simply a scaled version of the standard deviation $\sigma_t$. The
ratio $\lambda / \sigma_\epsilon$ acts as a normalizing constant that
accounts for the informativeness of observations.  The reason we
consider this scaling is that it highlights certain invariants across
problem instances. In Section~\ref{se:snr}, we will interpret the
value of $\rho_0$ as the relative volume of the trader's activity in
the marketplace. For the moment, it is sufficient to observe that if
the distribution $\phi_t$ is Gaussian, it is characterized by
$(\mu_t,\rho_t)$.

\section{Dynamic Programming Analysis}
\label{se:dp-analysis}

In this section, we develop abstract dynamic programming algorithms
for computing PBE and Gaussian PBE.  We also discuss structural
properties of associated value functions.  The dynamic programming
recursion relies on the computation of equilibria for single-stage
games, and we also discuss the existence of such equilibria.  The
algorithms of this section are not implementable, but their treatment
motivates the design of a practical algorithm that will be presented
in the next section.

\subsection{Stage-Wise Decomposition}\label{se:stage-wise}

The process of computing a PBE and the corresponding value functions
can be decomposed into a series of single-stage equilibrium problems
via a dynamic programming backward recursion. We begin by defining
some notation. For each $\pi_t$, $\psi_t$, and $u_t$, define a dynamic
programming operator $F^{(\psi_t,\hat{\pi}_t)}_{u_t}$ by
\[
\bigl(F^{(\psi_t, \hat{\pi}_t)}_{u_t} U\bigr)(x_{t-1},y_{t-1},\phi_{t-1})
\defeq \E^{(\psi_t, \hat{\pi}_t)}_{u_t}
\big[\left.\lambda (u_t + v_t) x_{t-1} 
  + U(x_t, y_t, \phi_t)\ \right|\  
x_{t-1}, y_{t-1}, \phi_{t-1}\big],
\]
for all functions $U$, where $x_t = x_{t-1} + u_t$, $y_t = y_{t-1} + v_t$, $v_t
= \psi_t(y_{t-1},\phi_{t-1})$, and $\phi_t$ results from the Bayesian
update \eqref{eq:bayes-update} given that the arbitrageur assumes the
trader trades $\hat{\pi}_t(x_{t-1},y_{t-1},\phi_{t-1})$ while the
trader actually trades $u_t$.  Similarly, for each $\pi_t$ and $v_t$,
define a dynamic programming operator $G_{v_t}^{\pi_t}$ by
\[
\bigl(G_{v_t}^{\pi_t} V\bigr)(y_{t-1},\phi_{t-1}) 
\defeq \E^{\pi_t}_{v_t}
\big[\left.\lambda (u_t + v_t) y_{t-1} + V(y_t,\phi_t)
  \ \right|\ y_{t-1}, \phi_{t-1}\big],
\]
for all functions $V$, where $y_t = y_{t-1} + v_t$, $u_t =
\pi_t(x_{t-1},y_{t-1},\phi_{t-1})$, $x_{t-1}$ is distributed
according to the belief $\phi_{t-1}$, and $\phi_t$ results from the
Bayesian update \eqref{eq:bayes-update} given that the arbitrageur
correctly assumes the trader trades $u_t$.

\begin{singlespace}
\begin{algorithm}[htbp]
\caption{PBE Solver}
\label{alg:pbe}
\begin{algorithmic}[1]
\STATE Initialize the terminal value functions $U^*_{T-1}$ and $V^*_{T-1}$ according to \eqref{eq:terminal-V}--\eqref{eq:terminal-U}\label{algline:pbe:terminal}
\FOR{$t=T-1,T-2,\ldots,1$} 
\STATE Compute $(\pi^*_t, \psi^*_t)$ such that for all $x_{t-1}$,
$y_{t-1}$, and $\phi_{t-1}$,\label{algline:pbe:pipsi}
\algeqstart
\begin{align*}
\pi^*_t(x_{t-1},y_{t-1},\phi_{t-1}) & \in 
\argmax_{u_t} \left(F^{(\psi^*_t, \pi^*_t)}_{u_t} U^*_t\right)(x_{t-1},y_{t-1},\phi_{t-1})
\\
\psi^*_t(y_{t-1},\phi_{t-1}) & \in \argmax_{v_t} \left(G_{v_t}^{\pi^*_t} V^*_t\right)(y_{t-1},\phi_{t-1})
\end{align*}
\algeqend
\STATE Compute the value functions at the previous time step by
setting, for all $x_{t-1}$, $y_{t-1}$, and $\phi_{t-1}$,\label{algline:pbe:UV}
\algeqstart
\vspace*{-\baselineskip}
\begin{align*}
U^*_{t-1}(x_{t-1},y_{t-1},\phi_{t-1}) 
& \leftarrow 
\left(F^{(\psi^*_t, \pi^*_t)}_{\pi_t^*} U^*_t\right)(x_{t-1},y_{t-1},\phi_{t-1})
\\
V^*_{t-1}(y_{t-1},\phi_{t-1}) 
& \leftarrow \left(G_{\psi_t^*}^{\pi^*_t} V^*_t\right)(y_{t-1},\phi_{t-1})
\end{align*}
\algeqend
\ENDFOR
\end{algorithmic}
\end{algorithm}
\end{singlespace}

Consider Algorithm~\ref{alg:pbe} for computing a PBE. In
Step~\ref{algline:pbe:terminal}, the algorithm begins by initializing
the terminal value functions $U^*_{T-1}$ and $V^*_{T-1}$. These terminal value functions have a simple closed form in equilibrium. This is because,
at time $T$, the trader must liquidate his position, hence
$\pi^*_T(x_{T-1},y_{T-1},\phi_{T-1})=-x_{T-1}$. Similarly, arbitrageur
must liquidate his position over times $T$ and $T+1$. In equilibrium,
he will do so optimally, thus his value function takes the form
\begin{equation}\label{eq:terminal-V}
  \begin{split}
    V^*_{T-1}(y_{T-1}, \phi_{T-1}) 
    & = \max_{v_T}\ 
    \E\left[\left.\lambda (- x_{T-1} + v_{T}) y_{T-1} 
        - \lambda (y_{T-1} + v_{T})^2
        \ \right|\ y_{T-1}, \phi_{T-1}\right]
    \\
    & = - \lambda \left( \mu_{T-1}
      + \tfrac{3}{4} y_{T-1}\right) y_{T-1},
  \end{split}
\end{equation}
where the optimizing decision is $\psi^*_T(y_{T-1},\phi_{T-1}) = 
-\tfrac{1}{2}y_{T-1}$. It is straightforward to derive the
corresponding expression of the trader's value function,
\begin{equation}\label{eq:terminal-U}
  \begin{split}
    U^*_{T-1}(x_{T-1},y_{T-1}, \phi_{T-1}) 
    & = \E\left[\left.\lambda \left(- x_{T-1} - \tfrac{1}{2} y_{T-1}\right)
        x_{T-1} 
          \ \right|\ x_{T-1}, y_{T-1}, \phi_{T-1}\right]
    \\
    & = - \lambda \left(x_{T-1} + \tfrac{1}{2} y_{T-1}\right) x_{T-1}.
  \end{split}
\end{equation}

At each time $t < T$, equilibrium policies must satisfy the
best-response conditions
\eqref{eq:trader-response}--\eqref{eq:arbitrageur-response}. Given the
value functions $U^*_t$ and $V^*_t$, these conditions decompose
recursively according to to Step~\ref{algline:pbe:pipsi}.  Given such
a pair $(\pi^*_t,\psi^*_t)$, the value functions $U^*_{t-1}$ and
$V^*_{t-1}$ for the prior time period are, in turn, computed in
Step~\ref{algline:pbe:UV}.
  
It is easy to see that, so long as
Step~\ref{algline:pbe:pipsi} is carried out successfully each time it
is invoked, the algorithm produces a PBE $(\pi^*, \phi^*)$ along with
value functions $U^*_t = U_t^{\pi^*, \psi^*}$ and $V^*_t =
V_t^{\psi^*,\pi^*}$.  However, the algorithm is not implementable.
For starters, the functions $\pi_t^*$, $\psi_t^*$, $U_{t-1}^*$, and
$V_{t-1}^*$, which must be computed and stored, have infinite domains.

\subsection{Linear Policies}
\label{se:lgpbe-recursion}

Consider the following class of policies:
\begin{definition}\label{def:linearpolicy}
A function $\pi_t$ is {\bf linear} if there are coefficients $a_{x,t}^{\rho_{t-1}}$, 
$a_{y,t}^{\rho_{t-1}}$ and $a_{\mu,t}^{\rho_{t-1}}$, which are functions of $\rho_{t-1}$, such that
\begin{equation}\label{eq:linear-pi}
\pi_t(x_{t-1},y_{t-1},\phi_{t-1}) = a_{x,t}^{\rho_{t-1}} x_{t-1} + a_{y,t}^{\rho_{t-1}} y_{t-1} + a_{\mu,t}^{\rho_{t-1}} \mu_{t-1},
\end{equation}
for all $x_{t-1}$, $y_{t-1}$, and $\phi_{t-1}$.  Similarly, function $\psi_t$ is {\bf linear} if there is a 
coefficients $b_{y,t}^{\rho_{t-1}}$ and $b_{\mu,t}^{\rho_{t-1}}$, 
which is a function of $\rho_{t-1}$, such that
\begin{equation}\label{eq:linear-psi}
\psi_t(y_{t-1},\phi_{t-1}) = b_{y,t}^{\rho_{t-1}} y_{t-1} + b_{\mu,t}^{\rho_{t-1}} \mu_{t-1},
\end{equation}
for all $y_{t-1}$ and $\phi_{t-1}$.  A policy is linear if the
component functions associated with times $1,\ldots,T-1$ are linear.
\end{definition}


By restricting attention to linear policies and Gaussian beliefs,
we can apply an algorithm similar to that presented in the previous
section to compute a Gaussian PBE.  In particular, consider Algorithm
\ref{alg:lgpbe}.  This algorithm aims to computes a single-stage
equilibrium that is linear.  Further, actions and values are only
computed and stored for elements of the domain for which $\phi_{t-1}$
is Gaussian.  This is only viable if the iterates $U_t^*$ and $V_t^*$,
which are computed only for Gaussian $\phi_t$, provide sufficient
information for subsequent computations.  This is indeed the case, as
a consequence of the following result.

\begin{theorem}\label{th:Gaussian}
  If the belief distribution $\phi_{t-1}$ at time is Gaussian, and the
  arbitrageur assumes that the trader's policy $\hat{\pi}_t$ is linear with
  $\hat{\pi}_t(x_{t-1},y_{t-1},\phi_{t-1}) = \hat{a}_{x,t}^{\rho_{t-1}}
  x_{t-1} + \hat{a}_{y,t}^{\rho_{t-1}} y_{t-1} + \hat{a}_{\mu,t}^{\rho_{t-1}}
  \mu_{t-1}$, then the belief distribution $\phi_t$ is also Gaussian. The mean
  $\mu_t$ is a linear function of $y_{t-1}$, $\mu_{t-1}$, and the observed
  price change $\Delta p_t$, with coefficients that are deterministic
  functions of the scaled variance $\rho_{t-1}$. The scaled variance $\rho_t$
  evolves according to
\begin{equation}\label{eq:rho-fwd-update}
\rho_t^2 = \left(1+\hat{a}^{\rho_{t-1}}_{x,t}\right)^2
\left(\frac{1}{\rho_{t-1}^2}
+ (\hat{a}^{\rho_{t-1}}_{x,t})^2
\right)^{-1}.
\end{equation}
In particular, $\rho_t$ is a deterministic function of $\rho_{t-1}$.
\end{theorem}
\noindent It follows from this result that if $\pi^*$ is linear
then, for Gaussian $\phi_{t-1}$, $F_{u_t}^{(\psi^*,\pi^*)} U^*_t$ only
depends on values of $U^*_t$ evaluated at Gaussian $\phi_t$.
Similarly, if $\pi^*$ is linear then, for Gaussian $\phi_{t-1}$,
$G_{v_t}^{\pi^*} V^*_t$ only depends on values of $V^*_t$ evaluated at
Gaussian $\phi_t$.  It also follows from this theorem that
Algorithm~\ref{alg:lgpbe}, which only computes actions and values for
Gaussian beliefs, results in a Gaussian PBE $(\pi^*,\psi^*)$.  We
should mention, though, that Algorithm~\ref{alg:lgpbe} is still not
implementable since the restricted domains of $U^*_t$ and $V^*_t$
remain infinite.

\begin{singlespace}
\begin{algorithm}[htbp]
\caption{Linear-Gaussian PBE Solver}
\label{alg:lgpbe}
\begin{algorithmic}[1]
\STATE Initialize the terminal value functions $U^*_{T-1}$ and $V^*_{T-1}$ according to \eqref{eq:terminal-V}--\eqref{eq:terminal-U}\label{algline:lgpbe:terminal}
\FOR{$t=T-1,T-2,\ldots,1$} 
\STATE Compute linear $(\pi^*_t, \psi^*_t)$ such that for all $x_{t-1}$,
$y_{t-1}$, and Gaussian $\phi_{t-1}$,\label{algline:lgpbe:pipsi}
\algeqstart
\begin{align*}
\pi^*_t(x_{t-1},y_{t-1},\phi_{t-1}) & \in 
\argmax_{u_t} \left(F^{(\psi^*_t, \pi^*_t)}_{u_t} U^*_t\right)(x_{t-1},y_{t-1},\phi_{t-1})
\\
\psi^*_t(y_{t-1},\phi_{t-1}) & \in \argmax_{v_t} \left(G_{v_t}^{\pi^*_t} V^*_t\right)(y_{t-1},\phi_{t-1})
\end{align*}
\algeqend
\STATE Compute the value functions at the previous time step by
setting, for all $x_{t-1}$, $y_{t-1}$, and Gaussian $\phi_{t-1}$,
\label{algline:lgpbe:UV}
\algeqstart
\begin{align*}
U^*_{t-1}(x_{t-1},y_{t-1},\phi_{t-1}) 
& \leftarrow 
\left(F^{(\psi^*_t, \pi^*_t)}_{\pi_t^*} U^*_t\right)(x_{t-1},y_{t-1},\phi_{t-1})
\\
V^*_{t-1}(y_{t-1},\phi_{t-1}) 
& \leftarrow \left(G_{\psi_t^*}^{\pi^*_t} V^*_t\right)(y_{t-1},\phi_{t-1})
\end{align*}
\algeqend
\ENDFOR
\end{algorithmic}
\end{algorithm}
\end{singlespace}

Motivated by these observations, for the remainder of the paper, we
will focus on computing equilibria of the following form:
\begin{definition}
  A pair of policies $(\pi^*,\psi^*)\in\Pi\times\Psi$ is a {\bf
    linear-Gaussian perfect Bayesian equilibrium} if it is a Gaussian
  PBE and each policy is linear.
\end{definition}

\subsection{Quadratic Value Functions}

Closely associated with linear policies are the following class of
value functions:
\begin{definition}
  A function $U_t$ is {\bf trader-quadratic-decomposable} (TQD) if
  there are coefficients $c_{xx,t}^{\rho_t}$, $c_{yy,t}^{\rho_t}$, $c_{\mu\mu,t}^{\rho_t}$,
  $c_{xy,t}^{\rho_t}$, $c_{x\mu,t}^{\rho_t}$, $c_{y\mu,t}^{\rho_t}$ and $c_{0,t}^{\rho_t}$, 
  which are functions of $\rho_t$, such that
\begin{equation}\label{eq:TQD}
\begin{split}
U_t(x_t,y_t,\phi_t) 
&
= -\lambda\bigg(\tfrac{1}{2} c^{\rho_t}_{xx,t} x_t^2 + \tfrac{1}{2} c^{\rho_t}_{yy,t} y_t^2
 + \tfrac{1}{2} c^{\rho_t}_{\mu\mu,t} \mu_t^2
+ c^{\rho_t}_{xy,t} x_t y_t + c^{\rho_t}_{x\mu,t} x_t \mu_t + c^{\rho_t}_{y\mu,t} y_t \mu_t 
- \frac{\sigma_\epsilon^2}{\lambda^2} c^{\rho_t}_{0,t}\bigg),
\end{split}
\end{equation}
for all $x_t$, $y_t$, and $\phi_t$.  A function $V_t$ as {\bf
  arbitrageur-quadratic-decomposable} (AQD) if there are coefficients
$d^{\rho_t}_{yy,t}$, $d^{\rho_t}_{\mu\mu,t}$, $d^{\rho_t}_{y\mu,t}$ 
and $d^{\rho_t}_{0,t}$, which are functions of $\rho_t$, such that
\begin{equation}\label{eq:AQD}
V_t(y_t,\phi_t) 
= 
- \lambda\left(\tfrac{1}{2} d^{\rho_t}_{yy,t} y_t^2 
+ \tfrac{1}{2} d^{\rho_t}_{\mu\mu,t} \mu_t^2
+ d^{\rho_t}_{y\mu,t} y_t \mu_t
-  \frac{\sigma_\epsilon^2}{\lambda^2} d^{\rho_t}_{0,t}\right),
\end{equation}
for all $y_t$ and $\phi_t$.
\end{definition}

In equilibrium, $U^*_{T-1}$ and $V^*_{T-1}$ are given by
Step~\ref{algline:lgpbe:terminal} of Algorithm~\ref{alg:lgpbe}, and
hence are TQD/AQD.  The following theorem captures how TQD and AQD
structure preserved in the dynamic programming recursion given linear
policies.
\begin{theorem}\label{th:tqd-aqd}
  If $U^*_t$ is TQD and $V^*_t$ is AQD, and
  Step~\ref{algline:lgpbe:pipsi} of Algorithm~\ref{alg:lgpbe} produces
  a linear pair $(\pi_t^*,\psi_t^*)$, then $U^*_{t-1}$ and
  $V^*_{t-1}$, defined by Step~\ref{algline:lgpbe:UV} of
  Algorithm~\ref{alg:lgpbe} are TQD and AQD, respectively.
\end{theorem}
\noindent Hence, each pair of value functions generated by
Algorithm~\ref{alg:lgpbe} is TQD/AQD.  A great benefit of this
property comes from the fact that, for a fixed value of $\rho_t$, each
associated value function can be encoded using just a few parameters.

\subsection{Simplified Conditions for Equilibrium}

Algorithm~\ref{alg:lgpbe} relies for each $t$ on existence of a pair
$(\pi^*_t, \psi^*_t)$ of linear functions that satisfy single-stage
equilibrium conditions. In general, this would require verifying that
each policy function is the Gaussian best response for all possible
states. The following theorem provides a much simpler set of
conditions. In Section~\ref{se:algorithm}, we will exploit these
conditions in order to compute equilibrium policies.

\begin{theorem}\label{th:sufficientcondition}
  Suppose that $U^*_t$ and $V^*_t$ and TQD/AQD value functions
  specified by \eqref{eq:TQD}--\eqref{eq:AQD}, and
  $(\pi^*_t,\psi^*_t)$ are linear policies specified by
  \eqref{eq:linear-pi}--\eqref{eq:linear-psi}.  Assume that, for all
  $\rho_{t-1}$, the policy coefficients satisfy the first order
  conditions
  \begin{gather}
    \begin{split}
      0 & =
      \big(\rho_t^2 c^{\rho_t}_{\mu\mu,t} + 2 \rho_t c^{\rho_t}_{x\mu,t} + c^{\rho_t}_{xx,t}\big)
      \big(a^{\rho_{t-1}}_{x,t}\big)^3 + \big(3c^{\rho_t}_{xx,t} + 3\rho_t
      c^{\rho_t}_{x\mu,t} - 1\big)\big(a^{\rho_{t-1}}_{x,t}\big)^2 \\
      & \quad + \big(3c^{\rho_t}_{xx,t} + \rho_t c^{\rho_t}_{x\mu,t} - 2\big)
      a^{\rho_{t-1}}_{x,t} + c^{\rho_t}_{xx,t} - 1, \\
    \end{split}\label{eq:simple-start}
    \\
    a^{\rho_{t-1}}_{y,t} 
    = 
    -\frac{\big(b^{\rho_{t-1}}_{y,t}+1\big) \big(c^{\rho_t}_{xy,t}+\alpha_t
      c^{\rho_t}_{y\mu,t}\big)}{c^{\rho_t}_{xx,t}+(\alpha_t +1) 
      c^{\rho_t}_{x \mu,t}+\alpha_t c^{\rho_t}_{\mu \mu,t}},\label{eq:recursion-ay} 
    \\
    a^{\rho_{t-1}}_{\mu,t} 
    =
    -\frac{a^{\rho_{t-1}}_{x,t} b^{\rho_{t-1}}_{\mu,t}
      \big(c^{\rho_t}_{xy,t}+\alpha_t c^{\rho_t}_{y \mu,t}\big)+\alpha_t 
      \big(c^{\rho_t}_{x \mu,t} + 
      \alpha_t c^{\rho_t}_{\mu\mu,t}\big)/\rho_{t-1}^2}
    {a^{\rho_{t-1}}_{x,t} \big(c^{\rho_t}_{xx,t}+(\alpha_t +1) 
      c^{\rho_t}_{x \mu,t}+\alpha_t c^{\rho_t}_{\mu \mu,t}\big)},\label{eq:recursion-amu} 
    \\
    b^{\rho_{t-1}}_{y,t} 
    =
    \frac{1-d^{\rho_t}_{y\mu,t}a^{\rho_{t-1}}_{y,t}}{d^{\rho_t}_{yy,t}} - 1,
    \qquad
    b^{\rho_{t-1}}_{\mu,t} 
    =
    -\frac{(1+a^{\rho_{t-1}}_{\mu,t} 
      + a^{\rho_{t-1}}_{x,t})d^{\rho_t}_{y\mu,t}}{d^{\rho_t}_{yy,t}},\label{eq:recursion-by-bmu}
  \end{gather}
  and the second order conditions
  \begin{equation}\label{eq:secondordercondition}
    c^{\rho_t}_{xx,t}+(\alpha_t +1) 
    c^{\rho_t}_{x \mu,t}+\alpha_t c^{\rho_t}_{\mu \mu,t} > 0,
    \qquad
    d^{\rho_t}_{yy,t} > 0,
  \end{equation}
  where the quantities $\alpha_t$ and $\rho_t$ satisfy
  \begin{equation}
    \label{eq:simple-end}
    \alpha_t = \frac{a^{\rho_{t-1}}_{x,t}\big(1 + a^{\rho_{t-1}}_{x,t}\big)}{1/\rho^2_{t-1} + \big(a^{\rho_{t-1}}_{x,t}\big)^2}, 
    \qquad
    \rho_t^2 = \left(1+a^{\rho_{t-1}}_{x,t}\right)^2
    \left(\frac{1}{\rho_{t-1}^2}
      + \big(a^{\rho_{t-1}}_{x,t}\big)^2
    \right)^{-1}.
  \end{equation}
  Then, $(\pi^*_t,\psi^*_t)$ satisfy the single-stage equilibrium conditions
  \begin{align*}
    \pi^*_t(x_{t-1},y_{t-1},\phi_{t-1}) & \in
    \argmax_{u_t} \left(F^{(\psi^*_t, \pi^*_t)}_{u_t}
      U^*_t\right)(x_{t-1},y_{t-1},\phi_{t-1}),
    \\
    \psi^*_t(y_{t-1},\phi_{t-1}) & \in 
    \argmax_{v_t} \left(G_{v_t}^{\pi^*_t}
      V^*_t\right)(y_{t-1},\phi_{t-1}),
  \end{align*}
  for all $x_{t-1}$, $y_{t-1}$, and Gaussian $\phi_{t-1}$.
\end{theorem}

Note that, while this theorem provides sufficient conditions for
linear policies satisfying equilibrium conditions, it does not
guarantee the existence or uniqueness of such policies. These remain
an open issues. However, we support the plausibility of existence
through the following result on Gaussian best responses to linear
policies. It asserts that, if $\psi_t$ and $\hat{\pi}_t$ are linear,
then there is a linear best-response $\pi_t$ for the trader in the
single-stage game.  Similarly, if $\pi_t$ is linear then there is a
linear best-response $\psi_t$ for the arbitrageur in the single-stage
game.
\begin{theorem}\label{th:linear-pi-psi}
  If $U_t$ is TQD, $\psi_t$ is linear, and $\hat{\pi}_t$ is
  linear, then there exists a linear $\pi_t$ such that
\[
\pi_t(x_{t-1}, y_{t-1},\phi_{t-1}) 
\in \argmax_{u_t}  
\left(F_{u_t}^{(\psi_t,\hat{\pi}_t)} U_t\right)(x_{t-1}, y_{t-1},\phi_{t-1}),
\]
for all $x_{t-1}$, $y_{t-1}$, and Gaussian $\phi_{t-1}$, so long as
the optimization problem is bounded.  Similarly, if $V_t$ is AQD and
$\pi_t$ is linear then there exists a linear $\psi_t$ such that
\[
\psi_t(y_{t-1},\phi_{t-1}) 
\in \argmax_{v_t} 
\left(G_{v_t}^{\pi_t} V_t\right)(y_{t-1},\phi_{t-1}),
\]
for all $y_{t-1}$ and Gaussian $\phi_{t-1}$, so long as the
optimization problem is bounded.
\end{theorem}
\noindent Based on these results, if the trader (arbitrageur) assumes
that the arbitrageur (trader) uses a linear policy then it
suffices for the trader (arbitrageur) to restrict himself to
linear policies.  Though not a proof of existence, this
observation that the set of linear policies is closed under the
operation of best response motivates an aim to compute
linear-Gaussian PBE.

\section{Algorithm}
\label{se:algorithm}

The previous section presented abstract algorithms and results that
lay the groundwork for the development of a practical algorithm which
we will present in this section.  We begin by discussing a
parsimonious representation of policies.

\subsection{Representation of Policies}\label{se:policy-rep}

Algorithm~\ref{alg:lgpbe} takes as input three values that
parameterize our model: $(\lambda,\sigma_\epsilon,T)$.  The algorithm
output can be encoded in terms of coefficients
$\{a^{\rho_{t-1}}_{x,t},
a^{\rho_{t-1}}_{y,t},a^{\rho_{t-1}}_{\mu,t},b^{\rho_{t-1}}_{y,t},
b^{\rho_{t-1}}_{\mu,t}\}$, for every $\rho_{t-1} > 0$ and each time
step\footnote{ Recall, from the discussion in
  Section~\ref{se:stage-wise}, that $a^{\rho_{t-1}}_{x,T} = -1$,
  $a^{\rho_{t-1}}_{y,T} = a^{\rho_{t-1}}_{\mu,T} = 0$,
  $b^{\rho_{t-1}}_{y,T} = -1/2$, $b^{\rho_{t-1}}_{\mu,T+1} =
  b^{\rho_{t-1}}_{\mu,T} = 0$, and $b^{\rho_{t-1}}_{y,T+1} = -1$, for
  all $\rho_{t-1}$.}  $t=1,\ldots,T-1$. These coefficients
parameterize linear-Gaussian PBE policies.  Note that the output
depends on $\lambda$ and $\sigma_\epsilon$ only through $\rho_t$.
Hence, given any $\lambda$ and $\sigma_\epsilon$ with the same
$\rho_t$, the algorithm obtains the same coefficients.  This means
that the algorithm need only be executed once to obtain solutions for
all choices of $\lambda$ and $\sigma_\epsilon$.

Now, for each $t$, the policy coefficients are deterministic functions
of $\rho_{t-1}$.  For a fixed value of $\rho_{t-1}$, the coefficients
can be stored as five numerical values.  However, it is not feasible
to simultaneously store coefficients associated with all possible
values of $\rho_{t-1}$. Fortunately, given a linear policy for the trader,
Theorem~\ref{th:Gaussian} establishes $\rho_t$ is a deterministic
function of $\rho_{t-1}$. Thus, the initial value $\rho_0$ determines all
subsequent values of $\rho_t$. It follows that, for a fixed value of
$\rho_0$, over the relevant portion of its domain, a linear-Gaussian
PBE can be encoded in terms of $5 (T-1)$ numerical values.  We will
design an algorithm that aims to compute these $5 (T-1)$ parameters,
which we will denote by
$\{a_{x,t},a_{y,t},a_{\mu,t},b_{y,t},b_{\mu,t}\}$, for
$t=1,\ldots,T-1$.  These parameters allow us to determine PBE actions
at all visited states, so long as the initial value of $\rho_0$ is
fixed.

\subsection{Searching for Equilibrium Variances}

The parameters $\{a_{x,t},a_{y,t},a_{\mu,t},b_{y,t},b_{y,t}\}$ characterize
linear-Gaussian PBE policies restricted to the sequence
$\rho_0,\ldots,\rho_{T-1}$ generated in the linear-Gaussian
PBE.  We do not know in advance what this sequence will be, and as
such, we seek simultaneously compute this sequence
alongside the policy parameters.

One way to proceed, reminiscent of the bisection method employed by
\citet{Kyle85} and \cite{Foster94} would be to conjecture a value for
$\rho_{T-1}$.  Given a candidate value $\hat{\rho}_{T-1}$, the
preceding values $\hat{\rho}_{T-2}, \ldots, \hat{\rho}_0$, along with
policy parameters for times $T-1,\ldots,1$, can be computed by
sequentially solving the equations
\eqref{eq:simple-start}--\eqref{eq:simple-end} for single-stage
equilibria. The resulting policies form a linear-Gaussian PBE,
restricted to the sequence $\hat{\rho}_0, \ldots, \hat{\rho}_{T-1}$
that they would generate if $\rho_0 = \hat{\rho}_0$. One can
then seek a value of $\hat{\rho}_{T-1}$ such that the resulting
$\hat{\rho}_0$ is indeed equal to $\rho_0$.  This can be accomplished,
for example, via bisection search.

The bisection method can be numerically unstable, however. This is
because, the belief update equation \eqref{eq:rho-fwd-update} is used to
sequentially compute the values
$\hat{\rho}_{T-2},\ldots,\hat{\rho}_{0}$ backwards in time. When the
target value of $\rho_0$ is very large, small changes in
$\hat{\rho}_{T-1}$ can result in very large changes in $\hat{\rho}_0$,
making it difficult to match the precisely value of $\rho_0$.

To avoid this numerical instability, consider
Algorithm~\ref{alg:lgpbe-search}. This algorithm maintains a guess
$\hat{\pi}$ of the equilibrium policy of the trader, and, along with
the initial value $\rho_0$, this is used to generate the sequence
$\hat{\rho}_1,\ldots,\hat{\rho}_{T-1}$ by applying the belief update
equation \eqref{eq:rho-fwd-update} {\em forward} in time. This
sequence of values is then used in the single-stage equilibrium
conditions to solve for policies $(\pi^*,\psi^*)$. A sequence of
values $\hat{\rho}_1,\ldots,\hat{\rho}_{T-1}$ is then computed forward
in time using the policy $\pi^*$. If this sequence matches the
sequence generated by the guess $\hat{\pi}$, then the algorithm has
converged. Otherwise, the algorithm is repeated with a new guess
policy that is a convex combination of $\hat{\pi}$ and $\pi^*$. Since this algorithm only ever applies the belief equation \eqref{eq:rho-fwd-update} forward in time, it does not suffer from the numerical instabilities of the bisection method.

Note that Step~\ref{alg:lgpbe-search:pipsi} of the algorithm treats
$\rho_{t-1}$ as a free variable that is solved alongside the policy parameters
$\{a_{x,t},a_{y,t},a_{\mu,t},b_{y,t},b_{\mu,t}\}$.  These variables are
computed by simultaneously solving the system of equations
\eqref{eq:simple-start}--\eqref{eq:simple-end} for single-stage
equilibrium. To be precise, $a_{x,t}$ is obtained by solving the cubic
polynomial equation \eqref{eq:simple-start} numerically.  Given a value for
$a_{x,t}$, the remaining parameters $\{a_{y,t},a_{\mu,t},b_{y,t},b_{\mu,t}\}$
are be obtained by solving the linear system of equations
\eqref{eq:recursion-ay}--\eqref{eq:recursion-by-bmu}, while $\rho_{t-1}$ is
obtained through \eqref{eq:simple-end} . It can then be verified that the
second order condition \eqref{eq:secondordercondition} holds.
Algorithm~\ref{alg:lgpbe-search} is implementable
and we use it in computational studies presented in the next section.

\begin{singlespace}
\begin{algorithm}[htbp]
\caption{Linear-Gaussian PBE Solver with Variance Search}
\label{alg:lgpbe-search}
\begin{algorithmic}[1]
\STATE Initialize $\hat{\pi}$ to an equipartitioning policy
\FOR{$k=1,2,\dots$}
\STATE Compute $\hat{\rho}_1,\ldots,\hat{\rho}_{T-1}$ according to the initial value $\rho_0$ and the policy $\hat{\pi}$ by \eqref{eq:rho-fwd-update}
\STATE Initialize the terminal value functions $U^*_{T-1}$ and $V^*_{T-1}$ according to \eqref{eq:terminal-V}--\eqref{eq:terminal-U}
\FOR{$t=T-1,T-2,\ldots,1$} 
\STATE Compute linear $(\pi^*_t, \psi^*_t)$ and $\rho_{t-1}$ solving the single-stage equilibrium conditions \eqref{eq:simple-start}--\eqref{eq:simple-end},
assuming that $\rho_t=\hat{\rho_t}$\label{alg:lgpbe-search:pipsi}
\STATE Compute the value functions $U^*_{t-1}$ and $V^*_{t-1}$ at the previous time step given $(\pi^*_t, \psi^*_t)$
\ENDFOR
\STATE Compute $\tilde{\rho}_1,\ldots,\tilde{\rho}_{T-1}$ according to the initial value $\rho_0$ and the policy $\pi^*$ by \eqref{eq:rho-fwd-update}
\IF{$\hat{\rho} = \tilde{\rho}$}
\RETURN
\ELSE
\STATE Set $\hat{\pi} \leftarrow \gamma_k\hat{\pi} + (1-\gamma_k) \pi^*$, where $\gamma_k\in [0,1)$ is a step-size
\ENDIF
\ENDFOR
\end{algorithmic}
\end{algorithm}
\end{singlespace}

\section{Computational Results}\label{se:computational}
\label{se:computation}

In this section, we present computational results generated using
Algorithm~\ref{alg:lgpbe-search}.  In Section~\ref{se:alt-policy}, we
introduce some alternative, intuitive policies which will serve as a basis of
comparison to the linear-Gaussian PBE policy. In Section~\ref{se:snr}, we
discuss the importance of the parameter $\rho_0\defeq \lambda \sigma_0 /
\sigma_\epsilon$ in the qualitative behavior of the Gaussian PBE policy and
interpret $\rho_0$ as a measure of the ``relative volume'' of the trader's
activity in the marketplace. In Section~\ref{se:rec}, we discuss the relative
performance of the policies from the perspective of the execution cost of the
trader. Here, we demonstrate experimentally that the Gaussian PBE policy can
offer substantial benefits. In Section~\ref{se:signal}, we examine the
signaling that occurs through price movements. Finally, in
Section~\ref{se:dyn-trade}, we highlight the fact that the PBE policy is
adaptive and dynamic, and seeks to exploit exogenous market fluctuations in
order to minimize execution costs.

\subsection{Alternative Policies}\label{se:alt-policy}

In order to understand the behavior of linear-Gaussian PBE policies,
we first define two alternative policies for the trader for the
purpose of comparison.  In the absence of an arbitrageur, it is
optimal for the trader to minimize execution costs by partitioning his
position into $T$ equally sized blocks and liquidating them
sequentially over the $T$ time periods, as established by
\citet{BertsimasLo98}. We refer to the resulting policy
$\pi^{\text{EQ}}$ as an {\em equipartitioning} policy. It is defined
by
\[
\pi^{\text{EQ}}_t(x_{t-1},y_{t-1},\phi_{t-1}) \defeq - \frac{1}{T - t + 1} x_{t-1},
\]
for all $t$, $x_{t-1}$, $y_{t-1}$, and $\phi_{t-1}$.

Alternatively, the trader may wish to liquidate his position in a way
so as to reveal as little information as possible to the
arbitrageur. Trading during the final two time periods $T-1$ and $T$
does not reveal information to the arbitrageur in a fashion that can
be exploited. This is because, as discussed in
Section~\ref{se:stage-wise}, the arbitrageur's optimal trades at time
$T$ and $T+1$ are $v_T=-y_{T-1}/2$ and $v_{T+1}=-y_T$, respectively,
and these are independent of any belief of the arbitrageur with
respect to the trader's position. Given that the trader is free to
trade over these two time periods without any information leakage, it
is natural to minimize execution cost by equipartitioning over these
two time periods. Hence, define the {\em minimum revelation} policy
$\pi^{\text{MR}}$ to be a policy that liquidates the trader's position
evenly across only the last two time periods. That is,
\[
\pi^{\text{MR}}_t(x_{t-1},y_{t-1},\phi_{t-1}) \defeq
\begin{cases}
0 & \text{if $t < T-1$,} \\
- \frac{1}{2} x_{t-1} & \text{if $t=T-1$,} \\
- x_{t-1} & \text{if $t=T$,} 
\end{cases}
\]
for all $t$, $x_{t-1}$, $y_{t-1}$, and $\phi_{t-1}$.

\subsection{Relative Volume}\label{se:snr}

Observed in Section~\ref{se:policy-rep}, linear-Gaussian
PBE policies are determined as a function of the composite parameter
$\rho_0\defeq \lambda \sigma_0 / \sigma_\epsilon$. In order to
interpret this parameter, consider the dynamics of price changes,
\[
\Delta p_{t} = \lambda (u_t + v_t) + \epsilon_t,\quad\epsilon_t\sim N(0,\sigma_\epsilon^2).
\]
Here, $\epsilon_t$ is interpreted as the exogenous, random component
of price changes. Alternatively, one can imagine the random component of price changes are arising from the price impact of ``noise traders''. Denote by $z_t$ the total order flow from noise traders at time $t$, and consider a model where
\[
\Delta p_{t} = \lambda (u_t + v_t + z_t),\quad z_t\sim N(0,\sigma_z^2).
\]
If $\sigma_\epsilon=\lambda \sigma_z$, these two models are
equivalent. In that case,
\[
\rho_0 \defeq \frac{\lambda \sigma_0}{\sigma_\epsilon} =
\frac{\sigma_0}{\sigma_z}.
\]
In other words, $\rho_0$ can be interpreted as the ratio of the
uncertainty of the total volume of the trader's activity to the per
period volume of noise trading. As such, we refer to $\rho_0$ as the
{\em relative volume}.

We shall see in the following sections that, qualitatively, the
performance and behavior of Gaussian PBE policies are determined by the
magnitude of $\rho_0$. In the high relative volume regime, when
$\rho_0$ is large, either the initial position uncertainty $\sigma_0$
is very large or the volatility $\sigma_z$ of the noise traders is
very small. In these cases, from the perspective of the arbitrageur,
the trader's activity contributes a significant informative signal
which can be decoded in the context of less significant exogenous
random noise. Hence, the trader's activity early in the time horizon
reveals significant information which can be exploited by the
arbitrageur. Thus, it may be better for the trader to defer his
liquidation until the end of the time horizon.

Alternatively, in the low relative regime, when $\rho_0$ is small, the
arbitrageur cannot effectively distinguish the activity of the trader
from the noise traders in the market. Hence, the trader is free to
distribute his trades across the time horizon so as to minimize market
impact, without fear of front-running by the arbitrageur.



\subsection{Policy Performance}\label{se:rec}

Consider a pair of policies $(\pi,\psi)$, and assume that the
arbitrageur begins with a position $y_0 = 0$ and an initial belief
$\phi_0=N(0, \sigma_0^2)$. Given an initial position $x_0$, the
trader's expected profit is $U^{\pi,\psi}_0(x_0, 0, \phi_0)$. One
might imagine, however, that the initial position $x_0$ represents one
of many different trials where the trader liquidates positions. It
makes sense for this distribution of $x_0$ over trials to be
consistent with the arbitrageurs belief $\phi_0$, since this belief
could be based on past trials.  Given this distribution, averaging
over trials results in expected profit
$\E[U_0^{\pi,\psi}(x_0,0,\phi_0) \ |\ \phi_0]$.  Alternatively, if the
trader liquidates his entire position immediately, the expected profit
becomes $\E[- \lambda x_0^2 \ |\ \phi_0] = - \lambda \sigma_0^2$.  We
define the {\em trader's normalized expected profit} $\bar{U}(\pi,\psi)$ to be
the ratio of these two quantities. When the trader's value function is TQD, this takes the form
\[
\bar{U}(\pi,\psi)
\defeq
\frac{\E\left[\left. U_0^{\pi,\psi}(x_0,0,\phi_0) 
\ \right|\ \phi_0\right]}{\lambda \sigma_0^2}
=
- \frac{1}{2} c^{\rho_0}_{xx,0} + \frac{1}{\rho^2_0} c^{\rho_0}_{0,0},
\]
where $c^{\rho_0}_{xx,0}$ and $c^{\rho_0}_{0,0}$ are the trader's appropriate
value function coefficients at time $t=0$.

Analogously, the {\em arbitrageur's normalized expected profit}
$\bar{V}(\pi,\psi)$ is defined to be the expected profit of the
arbitrageur normalized by the expected immediate liquidating cost of
the trader. When the arbitrageur's value function is AQD, this takes
the form
\[
\bar{V}(\pi,\psi) \defeq
\frac{\E\left[\left. V_0^{\psi,\pi}(x_0,0,\phi_0) 
\ \right|\ \phi_0\right]}{\lambda \sigma_0^2}
=
\frac{1}{\rho^2_0} d^{\rho_0}_{0,0},
\]

Now, let $(\pi^*,\psi^*)$ denote a linear-Gaussian PBE. Since the
corresponding value functions are TQD/AQD, the normalized expected profits
depend on the parameters $\{\sigma_0,\lambda,\sigma_\epsilon\}$ only
through the relative volume parameter $\rho_0 \defeq \lambda \sigma_0
/ \sigma_\epsilon$.

Similarly, given the equipartitioning policy $\pi^{\text{EQ}}$, define
$\psi^{\text{EQ}}$ to be the optimal response of the arbitrageur to
the trader's policy $\pi^{\text{EQ}}$. This best response policy can
be computed by solving the linear-quadratic control problem
corresponding to \eqref{eq:arbitrageur-response}, via dynamic
programming. The policy takes the form
\[
\psi^{\text{EQ}}_t(y_{t-1},\mu_{t-1}) 
= 
\begin{cases}
\frac{-1}{T+2-t} y_{t-1} 
- \frac{(T-t)(T-t+3)}{2 (T+1-t)(T+2-t)} \mu_{t-1} 
& \text{if $1 \leq t \leq T$,} \\
-y_T & \text{otherwise.}
\end{cases}
\]
Using a similar argument as above, it is easy to see that
$\bar{U}(\pi^{\text{EQ}},\psi^{\text{EQ}})$ and
$\bar{V}(\pi^{\text{EQ}},\psi^{\text{EQ}})$ are also functions of the
parameter $\rho_0$.

Finally, given the minimum revelation policy $\pi^{\text{MR}}$, define
$\psi^{\text{MR}}$ to be the optimal response of the arbitrageur to
the trader's policy $\pi^{\text{MR}}$. It can be shown that, when $y_0
= 0$ and $\mu_0 = 0$, the best response of the arbitrageur to the
minimum revelation policy is to do nothing--since no information is
revealed by the trader in a useful fashion, there is no opportunity to
front-run. Hence,
\[
\bar{U}(\pi^{\text{MR}},\psi^{\text{MR}})
=
\frac{\E\left[\left. - \frac{1}{2} \lambda x_0^2 - \frac{1}{4} \lambda x_0^2 
\ \right|\ \phi_0\right]}{\lambda \sigma_0^2} = -\frac{3}{4},
\qquad
\bar{V}(\pi^{\text{MR}},\psi^{\text{MR}}) = 0.
\]

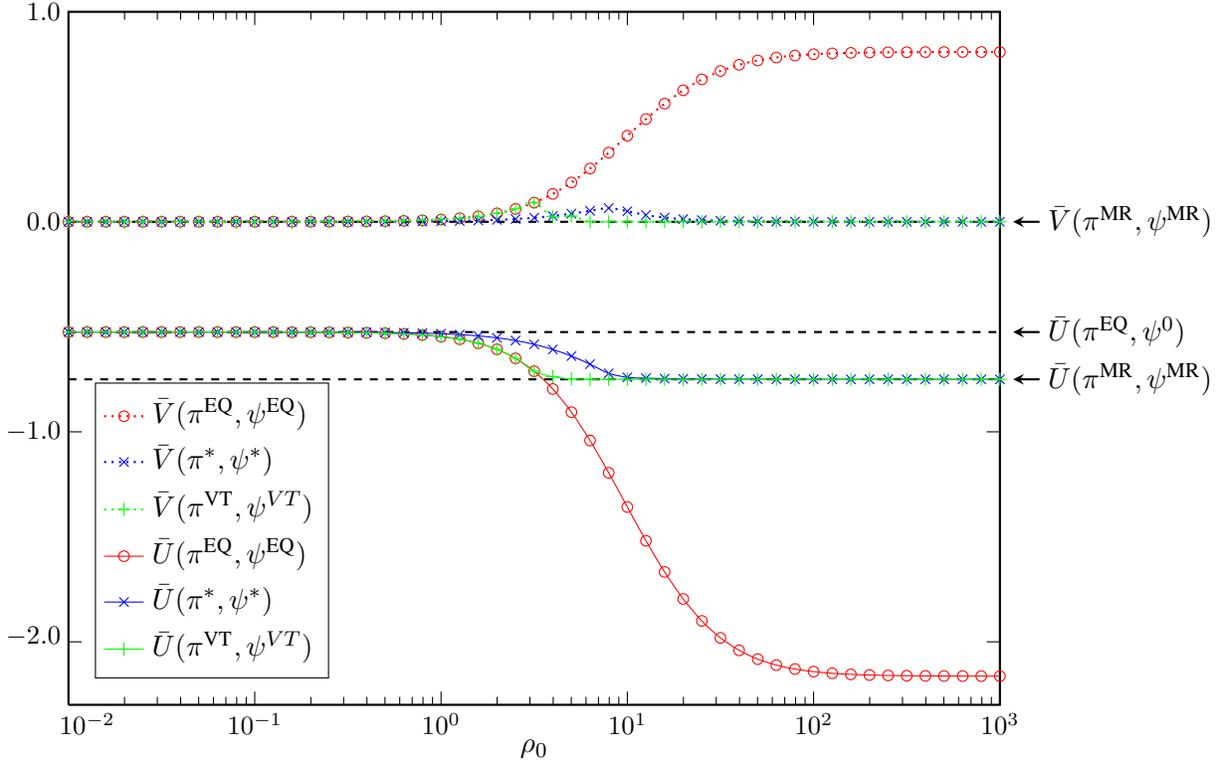
\begin{figure}[htb]
  \centering
  \hidefastcompile{
  \begin{tikzpicture}[x=0.975in,y=1.1in]
    \node[coordinate] at (-2,-2.3) (origin) {};
    \node[coordinate] at (3,0 |- origin) (x) {};
    \node[coordinate] at (0,1 -| origin) (y) {};
    \draw (0,-2 -| y) node[left,yshift=2pt] {\small $-2.0$};
    \foreach \y in {-1.0,0.0,1.0}
    \draw (0,\y -| y) node[left] {\small $\y$};
    \foreach \y in {-2.0,-1.0,0.0} {
      \draw[shift={(0,\y -| y)}] (5pt,0pt) -- (0pt,0pt);
      \draw[shift={(0,\y -| x)}] (-5pt,0pt) -- (0pt,0pt);
    }
    \draw (-2,0 |- x) node[below,xshift=7pt] {\small $10^{-2}$};
    \foreach \x in {-1, 0, 1, 2, 3}
    \draw (\x,0 |- x) node[below] {\small $10^{\x}$};
    \foreach \x in {-1, 0, 1, 2}
    {
      \draw[shift={(\x,0 |- x)}] (0pt,5pt) -- (0pt,0pt);
      \draw[shift={(\x,0 |- y)}] (0pt,-5pt) -- (0pt,0pt);
    }
    \foreach \x in {-1, 0, 1, 2, 3}
    {
      \foreach \z in {0.3010,0.4771,0.6021,0.6990,0.7782,0.8451,0.9031,0.9542}
      {
        \draw[shift={(\x,0 |- x)},shift={(-1,0)},shift={(\z,0)}]
        (0pt,3pt) -- (0pt,0pt);
        \draw[shift={(\x,0 |- y)},shift={(-1,0)},shift={(\z,0)}]
        (0pt,-3pt) -- (0pt,0pt);
      }
    }
    \draw[thick] (origin) to (x) node[midway,below=10pt] {$\rho_0$}
    to (y -| x) to (y) to (origin);
    \node[coordinate] at (0,-0.525) (eqna) {};
    \draw[thick,dashed] (eqna -| y) -- (eqna -| x)
    node[right=15pt] (REQNA) {$\bar{U}(\pi^{\text{EQ}},\psi^0)$};
    \draw[thick,->,shorten >=5pt] (REQNA.west) to (eqna -| x);
    \node[coordinate] at (0,-0.75) (Umr) {};
    \draw[thick,dashed] (Umr -| y) -- (Umr -| x)
    node[right=15pt] (RUMR) {$\bar{U}(\pi^{\text{MR}},\psi^{\text{MR}})$};
    \draw[thick,->,shorten >=5pt] (RUMR.west) to (RUMR -| x);
    \node[coordinate] at (0,0) (Vmr) {};
    \draw[thick,dashed] (Vmr -| y) -- (Vmr -| x)
    node[right=15pt] (RVMR) {$\bar{V}(\pi^{\text{MR}},\psi^{\text{MR}})$};
    \draw[thick,->,shorten >=5pt] (RVMR.west) to (RVMR -| x);
    \tikzstyle{Ucurve}=[smooth,solid,thin]
    \tikzstyle{Vcurve}=[smooth,dotted,thick]
    \tikzstyle{EQ}=[red,mark=o,mark size=2pt,mark options={thin,solid}]
    \tikzstyle{PBE}=[blue,mark=x,mark size=2.5pt,mark options={thin,solid}]
    \tikzstyle{VT}=[green,mark=+,mark size=2.5pt,mark options={thin,solid}]
    \draw[Ucurve,EQ]
    plot file {policy-performance-20-tr-eq.table};
    \draw[Ucurve,PBE]
    plot file {policy-performance-20-tr-pbe.table};
    \draw[Ucurve,VT]
    plot file {policy-performance-20-tr-heu.table};
    \draw[Vcurve,EQ]
    plot file {policy-performance-20-arb-eq.table};
    \draw[Vcurve,PBE]
    plot file {policy-performance-20-arb-pbe.table};
    \draw[Vcurve,VT]
    plot file {policy-performance-20-arb-heu.table};
    \matrix [draw,fill=white,column 2/.style={anchor=mid west},
    above right=10pt,ampersand replacement=\&] at (origin)
    {
      \draw[Vcurve,EQ,mark repeat=2,mark phase=2] plot coordinates 
      {([xshift=-7.5pt] 0,0) (0,0) ([xshift=7.5pt] 0,0) };
      \&
      \node {$\bar{V}(\pi^{\text{EQ}},\psi^{\text{EQ}})$};
      \\
      \draw[Vcurve,PBE,mark repeat=2,mark phase=2] plot coordinates 
      {([xshift=-7.5pt] 0,0) (0,0) ([xshift=7.5pt] 0,0) };
      \&
      \node {$\bar{V}(\pi^{*},\psi^{*})$};
      \\
      \draw[Vcurve,VT,mark repeat=2,mark phase=2] plot coordinates 
      {([xshift=-7.5pt] 0,0) (0,0) ([xshift=7.5pt] 0,0) };
      \&
      \node {$\bar{V}(\pi^{\text{VT}},\psi^{VT})$};
      \\
      \draw[Ucurve,EQ,mark repeat=2,mark phase=2] plot coordinates 
      {([xshift=-7.5pt] 0,0) (0,0) ([xshift=7.5pt] 0,0) };
      \&
      \node {$\bar{U}(\pi^{\text{EQ}},\psi^{\text{EQ}})$};
      \\
      \draw[Ucurve,PBE,mark repeat=2,mark phase=2] plot coordinates 
      {([xshift=-7.5pt] 0,0) (0,0) ([xshift=7.5pt] 0,0) };
      \&
      \node {$\bar{U}(\pi^{*},\psi^{*})$};
      \\
      \draw[Ucurve,VT,mark repeat=2,mark phase=2] plot coordinates 
      {([xshift=-7.5pt] 0,0) (0,0) ([xshift=7.5pt] 0,0) };
      \&
      \node {$\bar{U}(\pi^{\text{VT}},\psi^{VT})$};
      \\
    };
  \end{tikzpicture}
}
  \caption{The normalized expected profit of trading strategies for the time
    horizon $T=20$.\label{fig:policy-performance-20}}
\end{figure}

In Figure~\ref{fig:policy-performance-20}, the normalized expected profits of
various policies are plotted as functions of the relative volume
$\rho_0$, for a time horizon $T=20$. In all scenarios, as one might
expect, the trader's profit is negative while the arbitrageur's profit
is positive. In all cases, the trader's profit under the Gaussian PBE
policy dominates that under either the equipartitioning policy or the
minimum revelation policy. This difference is significant in moderate
to high relative volume regimes.

In the high relative volume regime, the equipartitioning policy fares
particularly badly from the perspective of the trader, performing up
to a factor of 2 worse than the Gaussian PBE policy. This effect
becomes more pronounced over longer time horizons. The minimum
revelation policy performs about as well as the PBE
policy. Asymptotically as $\rho_0 \uparrow\infty$, these policies
offer equivalent performance in the sense that
$\bar{U}(\pi^*,\psi^*)\uparrow
\bar{U}(\pi^{\text{MR}},\psi^{\text{MR}})=3/4$.

On the other hand, in the low relative volume regime, the
equipartitioning policy and the PBE policy perform comparably. Indeed,
define $\psi^0$ by $\psi^0_t \defeq 0$ for all $t$ (that is, no
trading by the arbitrageur). In the absence of an arbitrageur,
equipartitioning is the optimal policy for the trader, and backward
recursion can be used to show that
\[
\bar{U}(\pi^{\text{EQ}},\psi^0) = \frac{T+1}{2 T} \approx \frac{1}{2}.
\]
Asymptotically as $\rho_0 \downarrow 0$,
$\bar{U}(\pi^{\text{EQ}},\psi^{\text{EQ}}) \downarrow
\bar{U}(\pi^{\text{EQ}},\psi^0)$ and $\bar{U}(\pi^*,\psi^*) \downarrow
\bar{U}(\pi^{\text{EQ}},\psi^0)$. Thus, when the relative volume is
low, the effect of the arbitrageur becomes negligible when $\rho_0$ is
sufficiently small.

From the perspective of the arbitrageur in equilibrium,
$\bar{V}(\pi^*,\psi^*)\rightarrow 0$ as $\rho\rightarrow\pm\infty$. In the low
relative volume regime, the arbitrageur cannot distinguish the past activity
of the trader from noise, and hence is not able to profitably predict and
exploit the trader's future activity. In the high relative volume regime, as
we shall see in Section~\ref{se:dyn-trade}, the trader conceals his position
from the arbitrageur by deferring trading until the end of the horizon. Here,
as with the minimum revelation policy, the arbitrageur is not able to
profitably exploit the trader. Since the arbitrageur can choose not to trade
at each period, his best response to any trading strategy should lead to 
non-negative expected profit. In light of these observations, 
we can easily infer that in equilibrium the arbitrageur's profit curve should 
have at least one local maximum.

Both the equipartitioning and minimum revelation policies trade at a
constant rate, but over different, extremal time intervals: the
equipartitioning policy uses the entire time horizon, while the
minimum revelation policy uses only the last two time periods. A
fairer benchmark policy might consider optimizing the choice of time
interval. Define the {\em variable time} policy $\pi^{\text{VT}}$ as
follows: given the value $\rho_0$, select the $\tau$ such that trading
at a constant rate $u_t = -\frac{x_0}{\tau}$ over the last $\tau$ time
periods results in the highest expected profit for the trader,
assuming that the arbitrageur uses a best response policy. Define
$\psi^{\text{VT}}$ to be the best response of the arbitrageur to
$\pi^{\text{VT}}$. The variable time policy partially accounts for the
presence of an arbitrary, and the expected profit with the variable
time strategy will always be better that of equipartitioning or
minimum revelation. This is demonstrated by the
$\bar{U}(\pi^{\text{VT}},\psi^{\text{VT}})$ curve in
Figure~\ref{fig:policy-performance-20}. However, the trader still
fares better with an equilibrium policy, particularly in the
intermediate relative volume range, where the difference is close to
$20\%$.\footnote{In practice, improvements of as low as $0.01\%$ are
  considered significant.}

Examining Figure~\ref{fig:policy-performance-20}, it is clear that,
in equilibrium, the sum of the normalized profits of the trader and the
arbitrageur is negative, and the magnitude of sum is larger than the magnitude
of the loss incurred by the trader in the absence of the arbitrageur.  Define
the {\em spill-over}
to be the quantity
\[
\bar{U}(\pi^{\text{EQ}},\psi^0) - \left(\bar{U}(\pi^{*},\psi^{*}) + \bar{V}(\pi^{*},\psi^{*})\right).
\]
This is the difference between the normalized expected profit of the trader in
the absence of the arbitrageur, under the optimal equipartitioning policy, and
the combined normalized expected profits of the trader and arbitrageur in
equilibrium. The spill-over measures the benefit of the arbitrageur's presence
to the other participants of the system. Note that this benefit is positive,
and it is most significant in the high relative volume regime.

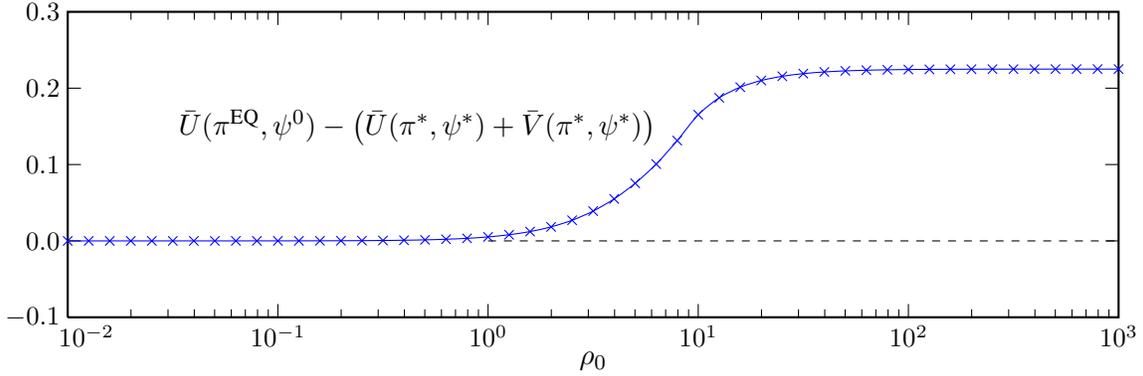
\begin{figure}[htb]
  \centering
  \hidefastcompile{
  \begin{tikzpicture}[x=1.1in,y=4in]
    \node[coordinate] at (-2,-0.1) (origin) {};
    \node[coordinate] at (3,0 |- origin) (x) {};
    \node[coordinate] at (0,0.3 -| origin) (y) {};
    \foreach \y in {-0.1,0.0,0.1,0.2,0.3}
    \draw (0,\y -| origin) node[left] {\small $\y$};
    \foreach \y in {-0.1,0.0,0.1,0.2,0.3} {
      \draw[shift={(0,\y -| y)}] (5pt,0pt) -- (0pt,0pt);
      \draw[shift={(0,\y -| x)}] (-5pt,0pt) -- (0pt,0pt);
    }
    \draw (-2,0 |- origin) node[below,xshift=7pt] {\small $10^{-2}$};
    \foreach \x in {-1, 0, 1, 2, 3} 
    \draw (\x,0 |- origin) node[below] {\small $10^{\x}$};
    \foreach \x in {-1, 0, 1, 2}
    {
      \draw[shift={(\x,0 |- x)}] (0pt,5pt) -- (0pt,0pt);
      \draw[shift={(\x,0 |- y)}] (0pt,-5pt) -- (0pt,0pt);
    }
    \foreach \x in {-1, 0, 1, 2, 3}
    {
      \foreach \z in {0.3010,0.4771,0.6021,0.6990,0.7782,0.8451,0.9031,0.9542}
      {
        \draw[shift={(\x,0 |- x)},shift={(-1,0)},shift={(\z,0)}]
        (0pt,3pt) -- (0pt,0pt);
        \draw[shift={(\x,0 |- y)},shift={(-1,0)},shift={(\z,0)}]
        (0pt,-3pt) -- (0pt,0pt);
      }
    }
    \draw[thick] (origin) to (x) node[midway,below=10pt] {$\rho_0$}
    to (y -| x) to (y) to (origin);
    \draw[dashed] (0,0) -- (0,0 -| x);
    \draw[smooth,solid,blue,mark=x,mark size=2.5pt,mark options={solid}] 
    plot file {policy-performance-20-efficiency.table};
    \node[yshift=5pt,left=5pt] at (0.90000,0.13138)
    {$\bar{U}(\pi^{\text{EQ}},\psi^0) - \left(\bar{U}(\pi^{*},\psi^{*}) 
        + \bar{V}(\pi^{*},\psi^{*})\right)$};
  \end{tikzpicture}
}
  \caption{The spill-over of the system for the time
    horizon $T=20$.\label{fig:policy-performance-20-eff}}
\end{figure}


In addition to the discussion of expected profits above, we can consider the variance of the trader's profits under different policies. Given a pair of policies $(\pi,\psi)$, define the {\em trader's normalized
  variance of profit} $\Var_U(\pi,\psi)$ as the variance under the
policies $(\pi,\psi)$ relative to the variance of immediate
liquidation. In other words,
\[
\Var_U(\pi,\psi) 
=
\frac{
  \Var^{\pi,\psi}\left(\left.\sum_{\tau=0}^{T-1} \Delta p_{\tau+1} x_\tau\ 
    \right|\ \phi_0\right)}
{\Var\left(\left. -\lambda x_0^2 + \epsilon_1 x_0\ 
      \right|\ \phi_0\right)}
=
\frac{
  \Var^{\pi,\psi}\left(\left.\sum_{\tau=0}^{T-1} \Delta p_{\tau+1} x_\tau\ 
    \right|\ \phi_0\right)}
{2 \lambda^2 \sigma_0^4 + \sigma_\epsilon^2
\sigma_0^2},
\]
where, as before, the expectations are taken assuming the policies
$(\pi,\psi)$ are used, $y_0 = \mu_0 = 0$, and
$x_0\sim\phi_0=N(0,\sigma_0^2)$. Similarly, it is possible to see
that, for a pair of linear policies $(\pi,\psi)$, the trader's
normalized variance of profit depends on the model parameters
$\{\sigma_0,\lambda,\sigma_\epsilon\}$ only through $\rho_0$.

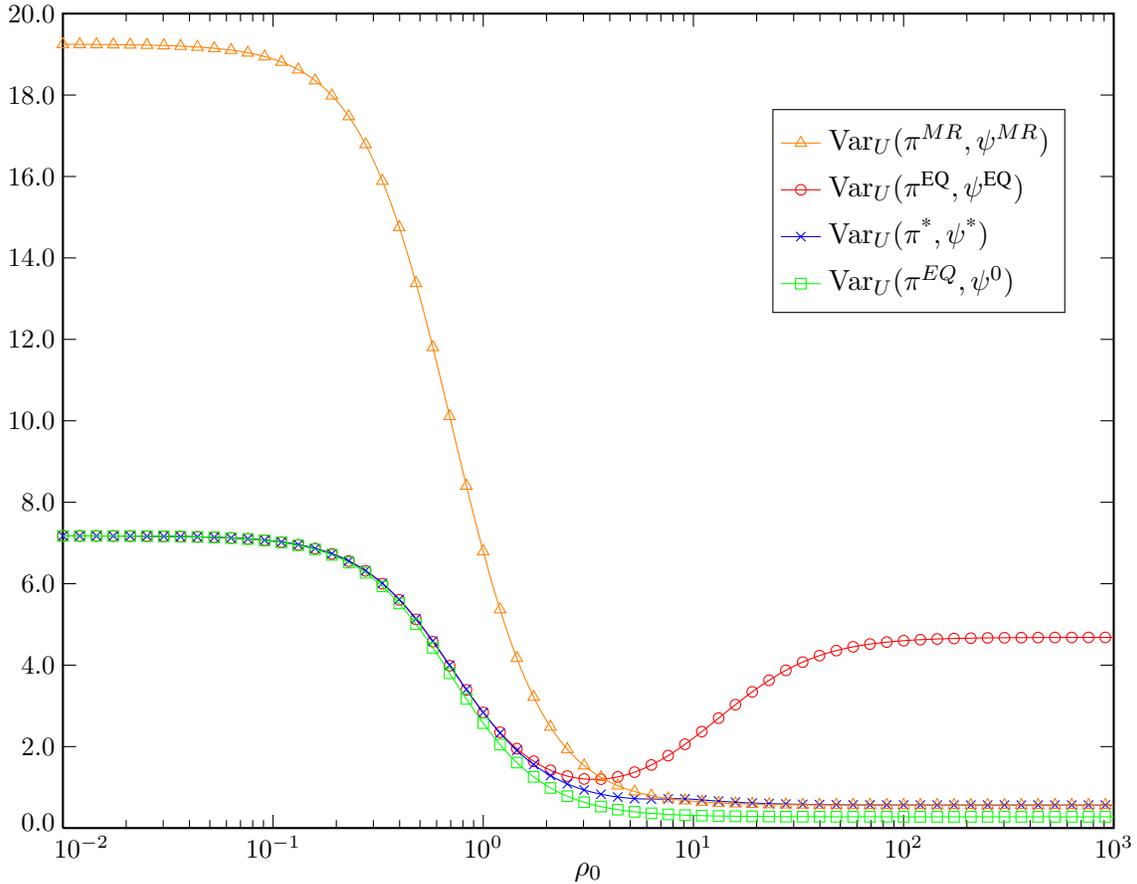
\begin{figure}[htb]
  \centering
  \hidefastcompile{
  \begin{tikzpicture}[x=1.1in,y=0.2133in]
    \node[coordinate] at (-2,0) (origin) {};
    \node[coordinate] at (3,0 |- origin) (x) {};
    \node[coordinate] at (0,20 -| origin) (y) {};
    \draw (0,0 -| y) node[left,yshift=2pt] {\small $0.0$};
    \foreach \y in {2.0,4.0,6.0,8.0,10.0,12.0,14.0,16.0,18.0,20.0} 
    \draw (0,\y -| y) node[left] {\small $\y$};
    \foreach \y in {2.0,4.0,6.0,8.0,10.0,12.0,14.0,16.0,18.0} {
      \draw[shift={(0,\y -| y)}] (5pt,0pt) -- (0pt,0pt);
      \draw[shift={(0,\y -| x)}] (-5pt,0pt) -- (0pt,0pt);
    }
    \draw (-2,0 |- x) node[below,xshift=7pt] {\small $10^{-2}$};
    \foreach \x in {-1, 0, 1, 2, 3} 
    \draw (\x,0 |- x) node[below] {\small $10^{\x}$};
    \foreach \x in {-1, 0, 1, 2} {
      \draw[shift={(\x,0 |- x)}] (0pt,5pt) -- (0pt,0pt);
      \draw[shift={(\x,0 |- y)}] (0pt,-5pt) -- (0pt,0pt);
    }
    \foreach \x in {-1, 0, 1, 2, 3}
    {
      \foreach \z in {0.3010,0.4771,0.6021,0.6990,0.7782,0.8451,0.9031,0.9542}
      {
        \draw[shift={(\x,0 |- x)},shift={(-1,0)},shift={(\z,0)}]
        (0pt,3pt) -- (0pt,0pt);
        \draw[shift={(\x,0 |- y)},shift={(-1,0)},shift={(\z,0)}]
        (0pt,-3pt) -- (0pt,0pt);
      }
    }
    \draw[thick] (origin) to (x) node[midway,below=10pt] {$\rho_0$}
    to (y -| x) to (y) to (origin)
    node[midway] (ylabel) {};
    \tikzstyle{Ucurve}=[smooth,solid,thin,mark repeat=4]
    \tikzstyle{EQ}=[red,mark=o,mark size=2pt,mark options={thin,solid}]
    \tikzstyle{EQ0}=[green,mark=square,mark size=2pt,mark options={thin,solid}]
    \tikzstyle{PBE}=[blue,mark=x,mark size=2.5pt,mark options={thin,solid}]
    \tikzstyle{MR}=[orange,mark=triangle,mark size=2.5pt,
    mark options={thin,solid}]
    \draw[Ucurve,EQ] plot file {normalized-variance-20-tr-eq.table};
    \draw[Ucurve,PBE] plot file {normalized-variance-20-tr-pbe.table};
    \draw[Ucurve,EQ0] plot file {normalized-variance-20-tr-eq0.table};
    \draw[Ucurve,MR] plot file {normalized-variance-20-tr-mr.table};
    \matrix [draw,fill=white,column 2/.style={anchor=mid west},
    above right=10pt,ampersand replacement=\&] at (1.25,12)
    {
      \draw[Ucurve,MR,mark repeat=2,mark phase=2] plot coordinates 
      {([xshift=-7.5pt] 0,0) (0,0) ([xshift=7.5pt] 0,0) };
      \&
      \node {$\Var_U(\pi^{MR},\psi^{MR})$};
      \\
      \draw[Ucurve,EQ,mark repeat=2,mark phase=2] plot coordinates 
      {([xshift=-7.5pt] 0,0) (0,0) ([xshift=7.5pt] 0,0) };
      \&
      \node {$\Var_U(\pi^{\text{EQ}},\psi^{\text{EQ}})$};
      \\
      \draw[Ucurve,PBE,mark repeat=2,mark phase=2] plot coordinates 
      {([xshift=-7.5pt] 0,0) (0,0) ([xshift=7.5pt] 0,0) };
      \&
      \node {$\Var_U(\pi^{\text{*}},\psi^{\text{*}})$};
      \\
      \draw[Ucurve,EQ0,mark repeat=2,mark phase=2] plot coordinates 
      {([xshift=-7.5pt] 0,0) (0,0) ([xshift=7.5pt] 0,0) };
      \&
      \node {$\Var_U(\pi^{EQ},\psi^{0})$};
      \\
    };
  \end{tikzpicture}  
}
  \caption{The trader's normalized variance of profit for the time
    horizon $T=20$. \label{fig:normalized-variance-20}}
\end{figure}

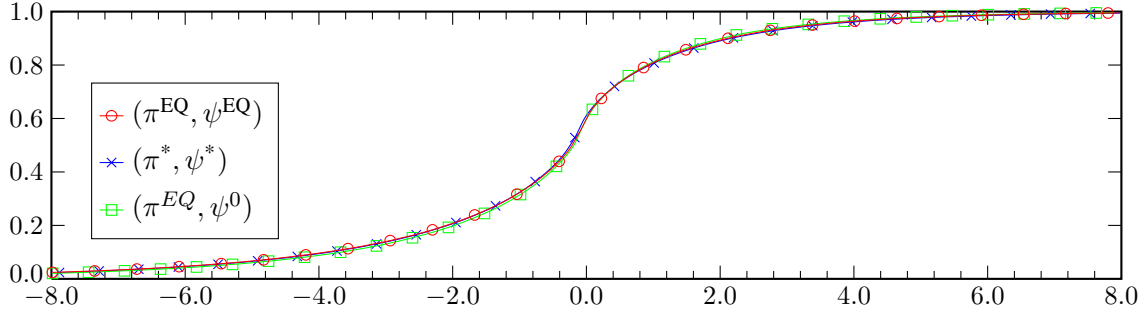
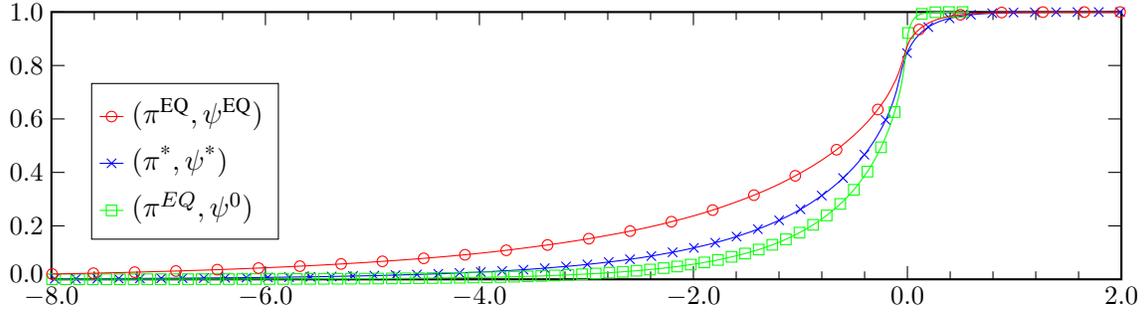
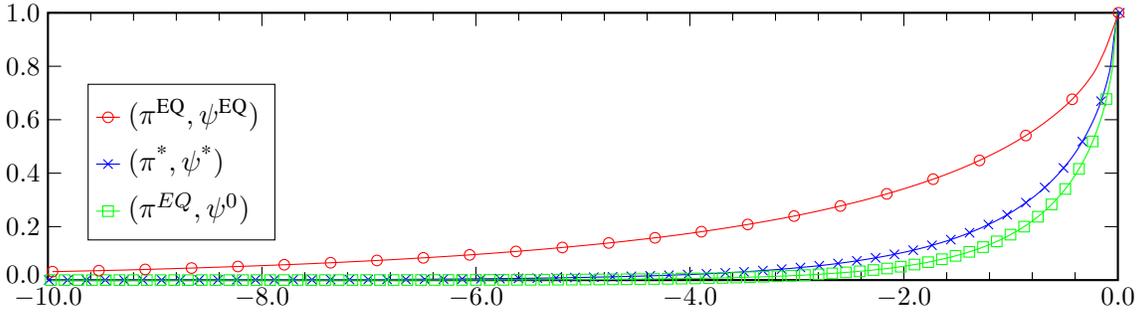
\begin{figure}[htb]
  \centering
  \subfigure[The low relative volume regime, $\rho_0 = 1$.]{
  \begin{tikzpicture}[x=0.35in,y=1.4in]
    \node[coordinate] at (-8,0) (origin) {};
    \node[coordinate] at (8,0 |- origin) (x) {};
    \node[coordinate] at (0,1 -| origin) (y) {};
    \draw (0,0 -| y) node[left,yshift=2pt] {\small $0.0$};
    \foreach \y in {0.2,0.4,0.6,0.8,1.0} 
    \draw (0,\y -| y) node[left] {\small $\y$};
    \foreach \y in {0.2,0.4,0.6,0.8,1.0} {
      \draw[shift={(0,\y -| y)}] (5pt,0pt) -- (0pt,0pt);
      \draw[shift={(0,\y -| x)}] (-5pt,0pt) -- (0pt,0pt);
    }
    \draw (-8,0 |- x) node[below,xshift=0pt] {\small $-8.0$};
    \foreach \x in {-6.0,-4.0,-2.0,0.0,2.0,4.0,6.0,8.0} 
    \draw (\x,0 |- x) node[below] {\small $\x$};
    \foreach \x in {-6.0,-4.0,-2.0,0.0,2.0,4.0,6.0} {
      \draw[shift={(\x,0 |- x)}] (0pt,5pt) -- (0pt,0pt);
      \draw[shift={(\x,0 |- y)}] (0pt,-5pt) -- (0pt,0pt);
    }
    \foreach \x in {-8.0,-6.0,-4.0,-2.0,0.0,2.0,4.0,6.0}
    {
      \foreach \z in {0.4,0.8,1.2,1.6}
      {
        \draw[shift={(\x,0 |- x)},shift={(0,0)},shift={(\z,0)}]
        (0pt,3pt) -- (0pt,0pt);
        \draw[shift={(\x,0 |- y)},shift={(0,0)},shift={(\z,0)}]
        (0pt,-3pt) -- (0pt,0pt);
      }
    }
    \draw[thick] (origin) to (x) node[midway,below=10pt] {}
    to (y -| x) to (y) to (origin)
    node[midway] (ylabel) {};
    \tikzstyle{Ucurve}=[smooth,solid,thin,mark repeat=4]
    \tikzstyle{Vcurve}=[smooth,dotted,thick,mark repeat=4]
    \tikzstyle{EQ}=[red,mark=o,mark size=2pt,mark options={thin,solid}]
    \tikzstyle{EQ0}=[green,mark=square,mark size=2pt,mark options={thin,solid}]
    \tikzstyle{PBE}=[blue,mark=x,mark size=2.5pt,mark options={thin,solid}]
    \draw[Ucurve,PBE] plot file {distribution-norm-profit-L-PBE-tr.table};
    \draw[Ucurve,EQ0] plot file {distribution-norm-profit-L-EQAlone-tr.table};
    \draw[Ucurve,EQ] plot file {distribution-norm-profit-L-EQEQ-tr.table};
    \matrix [draw,fill=white,column 2/.style={anchor=mid west},
    above right=15pt,ampersand replacement=\&] at (origin)
    {
      \draw[Ucurve,EQ,mark repeat=2,mark phase=2] plot coordinates 
      {([xshift=-4.5pt] 0,0) (0,0) ([xshift=4.5pt] 0,0) }; 
      \& 
      \node{$(\pi^{\text{EQ}},\psi^{\text{EQ}})$}; 
      \\
      \draw[Ucurve,PBE,mark repeat=2,mark phase=2] plot coordinates 
      {([xshift=-4.5pt] 0,0) (0,0) ([xshift=4.5pt] 0,0) }; 
      \& 
      \node{$(\pi^{\text{*}},\psi^{\text{*}})$}; 
      \\
      \draw[Ucurve,EQ0,mark repeat=2,mark phase=2] plot coordinates 
      {([xshift=-4.5pt] 0,0) (0,0) ([xshift=4.5pt] 0,0) };
      \& 
      \node{$(\pi^{EQ},\psi^{0})$}; 
      \\
    };
  \end{tikzpicture}
  }
  \subfigure[The moderate relative volume regime, $\rho_0 = 10$.]{
  \begin{tikzpicture}[x=0.56in,y=1.4in]
    \node[coordinate] at (-8,0) (origin) {};
    \node[coordinate] at (2,0 |- origin) (x) {};
    \node[coordinate] at (0,1 -| origin) (y) {};
    \draw (0,0 -| y) node[left,yshift=2pt] {\small $0.0$};
    \foreach \y in {0.2,0.4,0.6,0.8,1.0} 
    \draw (0,\y -| y) node[left] {\small $\y$};
    \foreach \y in {0.2,0.4,0.6,0.8,1.0} {
      \draw[shift={(0,\y -| y)}] (5pt,0pt) -- (0pt,0pt);
      \draw[shift={(0,\y -| x)}] (-5pt,0pt) -- (0pt,0pt);
    }
    \draw (-8,0 |- x) node[below,xshift=0pt] {\small $-8.0$};
    \foreach \x in {-6.0,-4.0,-2.0,0.0,2.0}
    \draw (\x,0 |- x) node[below] {\small $\x$};
    \foreach \x in {-6.0,-4.0,-2.0,0.0} {
      \draw[shift={(\x,0 |- x)}] (0pt,5pt) -- (0pt,0pt);
      \draw[shift={(\x,0 |- y)}] (0pt,-5pt) -- (0pt,0pt);
    }
    \foreach \x in {-8.0,-6.0,-4.0,-2.0,0.0}
    {
      \foreach \z in {0.4,0.8,1.2,1.6}
      {
        \draw[shift={(\x,0 |- x)},shift={(0,0)},shift={(\z,0)}]
        (0pt,3pt) -- (0pt,0pt);
        \draw[shift={(\x,0 |- y)},shift={(0,0)},shift={(\z,0)}]
        (0pt,-3pt) -- (0pt,0pt);
      }
    }
    \draw[thick] (origin) to (x) node[midway,below=10pt] {}
    to (y -| x) to (y) to (origin)
    node[midway] (ylabel) {};
    \tikzstyle{Ucurve}=[smooth,solid,thin,mark repeat=4]
    \tikzstyle{Vcurve}=[smooth,dotted,thick,mark repeat=4]
    \tikzstyle{EQ}=[red,mark=o,mark size=2pt,mark options={thin,solid}]
    \tikzstyle{EQ0}=[green,mark=square,mark size=2pt,mark options={thin,solid}]
    \tikzstyle{PBE}=[blue,mark=x,mark size=2.5pt,mark options={thin,solid}]
    \draw[Ucurve,PBE] plot file {distribution-norm-profit-M-PBE-tr.table};
    \draw[Ucurve,EQ0] plot file {distribution-norm-profit-M-EQAlone-tr.table};
    \draw[Ucurve,EQ] plot file {distribution-norm-profit-M-EQEQ-tr.table};
    \matrix [draw,fill=white,column 2/.style={anchor=mid west},
    above right=15pt,ampersand replacement=\&] at (origin)
    {
      \draw[Ucurve,EQ,mark repeat=2,mark phase=2] plot coordinates 
      {([xshift=-4.5pt] 0,0) (0,0) ([xshift=4.5pt] 0,0) }; 
      \& 
      \node{$(\pi^{\text{EQ}},\psi^{\text{EQ}})$}; 
      \\
      \draw[Ucurve,PBE,mark repeat=2,mark phase=2] plot coordinates 
      {([xshift=-4.5pt] 0,0) (0,0) ([xshift=4.5pt] 0,0) }; 
      \& 
      \node{$(\pi^{\text{*}},\psi^{\text{*}})$}; 
      \\
      \draw[Ucurve,EQ0,mark repeat=2,mark phase=2] plot coordinates 
      {([xshift=-4.5pt] 0,0) (0,0) ([xshift=4.5pt] 0,0) };
      \& 
      \node{$(\pi^{EQ},\psi^{0})$}; 
      \\
    };
  \end{tikzpicture}
  }
  \subfigure[The high relative volume regime, $\rho_0 = 100$.]{
  \begin{tikzpicture}[x=0.56in,y=1.4in]
    \node[coordinate] at (-10,0) (origin) {};
    \node[coordinate] at (0,0 |- origin) (x) {};
    \node[coordinate] at (0,1 -| origin) (y) {};
    \draw (0,0 -| y) node[left,yshift=2pt] {\small $0.0$};
    \foreach \y in {0.2,0.4,0.6,0.8,1.0} 
    \draw (0,\y -| y) node[left] {\small $\y$};
    \foreach \y in {0.2,0.4,0.6,0.8,1.0} {
      \draw[shift={(0,\y -| y)}] (5pt,0pt) -- (0pt,0pt);
      \draw[shift={(0,\y -| x)}] (-5pt,0pt) -- (0pt,0pt);
    }
    \draw (-10,0 |- x) node[below,xshift=0pt] {\small $-10.0$};
    \foreach \x in {-8.0,-6.0,-4.0,-2.0,0.0}
    \draw (\x,0 |- x) node[below] {\small $\x$};
    \foreach \x in {-8.0,-6.0,-4.0,-2.0} {
      \draw[shift={(\x,0 |- x)}] (0pt,5pt) -- (0pt,0pt);
      \draw[shift={(\x,0 |- y)}] (0pt,-5pt) -- (0pt,0pt);
    }
    \foreach \x in {-10.0,-8.0,-6.0,-4.0,-2.0}
    {
      \foreach \z in {0.4,0.8,1.2,1.6}
      {
        \draw[shift={(\x,0 |- x)},shift={(0,0)},shift={(\z,0)}]
        (0pt,3pt) -- (0pt,0pt);
        \draw[shift={(\x,0 |- y)},shift={(0,0)},shift={(\z,0)}]
        (0pt,-3pt) -- (0pt,0pt);
      }
    }
    \draw[thick] (origin) to (x) node[midway,below=10pt] {}
    to (y -| x) to (y) to (origin)
    node[midway] (ylabel) {};
    \tikzstyle{Ucurve}=[smooth,solid,thin,mark repeat=4]
    \tikzstyle{Vcurve}=[smooth,dotted,thick,mark repeat=4]
    \tikzstyle{EQ}=[red,mark=o,mark size=2pt,mark options={thin,solid}]
    \tikzstyle{EQ0}=[green,mark=square,mark size=2pt,mark options={thin,solid}]
    \tikzstyle{PBE}=[blue,mark=x,mark size=2.5pt,mark options={thin,solid}]
    \draw[Ucurve,PBE] plot file {distribution-norm-profit-H-PBE-tr.table};
    \draw[Ucurve,EQ0] plot file {distribution-norm-profit-H-EQAlone-tr.table};
    \draw[Ucurve,EQ] plot file {distribution-norm-profit-H-EQEQ-tr.table};
    \matrix [draw,fill=white,column 2/.style={anchor=mid west},
    above right=15pt,ampersand replacement=\&] at (origin)
    {
      \draw[Ucurve,EQ,mark repeat=2,mark phase=2] plot coordinates 
      {([xshift=-4.5pt] 0,0) (0,0) ([xshift=4.5pt] 0,0) }; 
      \& 
      \node{$(\pi^{\text{EQ}},\psi^{\text{EQ}})$}; 
      \\
      \draw[Ucurve,PBE,mark repeat=2,mark phase=2] plot coordinates 
      {([xshift=-4.5pt] 0,0) (0,0) ([xshift=4.5pt] 0,0) }; 
      \& 
      \node{$(\pi^{\text{*}},\psi^{\text{*}})$}; 
      \\
      \draw[Ucurve,EQ0,mark repeat=2,mark phase=2] plot coordinates 
      {([xshift=-4.5pt] 0,0) (0,0) ([xshift=4.5pt] 0,0) };
      \& 
      \node{$(\pi^{EQ},\psi^{0})$}; 
      \\
    };
  \end{tikzpicture}
  }
  \medskip

  \caption{The cumulative distribution of trader's normalized profit for 
    the time horizon $T=20$. \label{fig:profit-distribution-20}}
\end{figure}

In Figure~\ref{fig:normalized-variance-20}, the trader's normalized variance
of profit is plotted under the different policies. The lowest variance occurs
when the trader equipartitions and there is no arbitrageur, this is the curve
$\Var_U(\pi^{\text{EQ}},\psi^0)$. When the arbitrageur is present, however,
the variance in equilibrium $\Var_U(\pi^*,\psi^*)$ is less than either when the
trader equipartitions (i.e., the curve
$\Var_U(\pi^{\text{EQ}},\psi^{\text{EQ}})$) or employs the minimum revelation
policy (i.e., the curve $\Var_U(\pi^{\text{MR}},\psi^{\text{MR}})$).
Figure~\ref{fig:profit-distribution-20} shows the entire cumulative
distribution function of the trader's normalized profit under various relative
volume regimes. Given the presence of the arbitrageur, the equilibrium policy
has second-order dominance over equipartitioning in all relative volume
regimes.

\subsection{Signaling}\label{se:signal}

An important aspect of the linear-Gaussian PBE policy is that it accounts for
information conveyed through price movements. In order to understand
this feature, define the {\em relative uncertainty} to be the standard
deviation of the arbitrageur's belief about the trader's position at time
$t$, relative to that of the belief at time $0$; i.e., the ratio
$\sigma_t / \sigma_0$. By considering the evolution of relative
uncertainty over time for the Gaussian PBE policy versus the equipartitioning
and minimum revelation policies, we can study the comparative
signaling behavior.

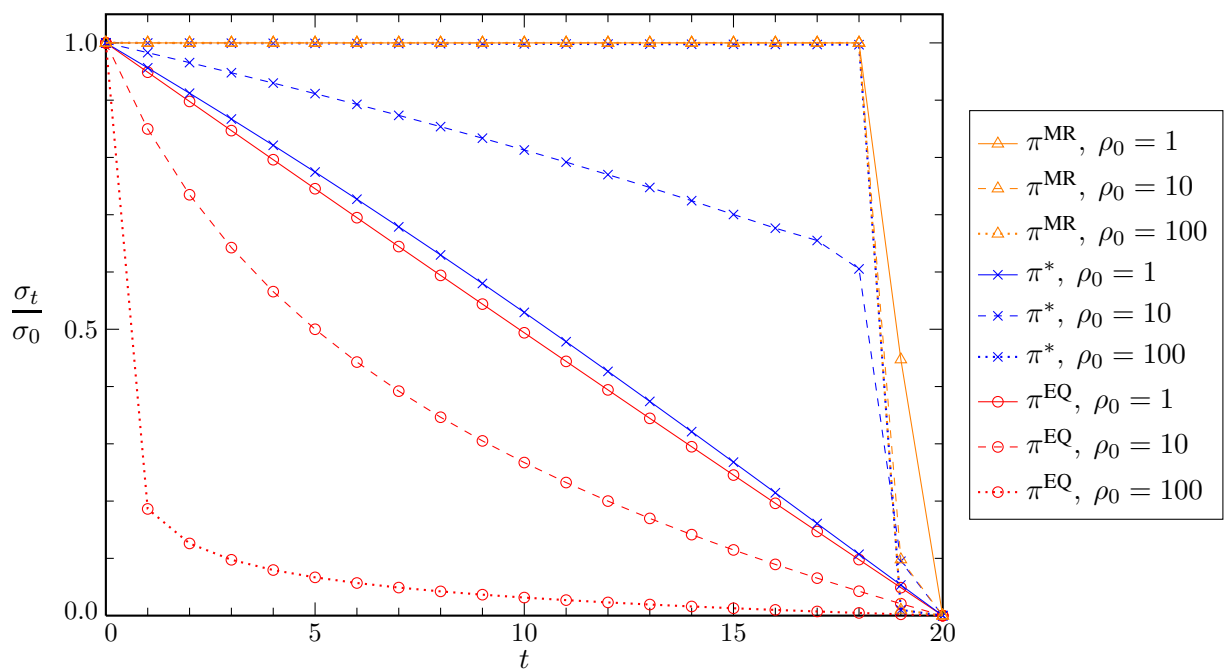
\begin{figure}[htb]
  \centering
  \hidefastcompile{
  \begin{tikzpicture}[x=0.219in,y=3.0in]
    \node[coordinate] at (0,0) (origin) {};
    \node[coordinate] at (20,0 |- origin) (x) {};
    \node[coordinate] at (0,1.05 -| origin) (y) {};
    \draw (0,0 -| origin) node[left,yshift=2pt] {\small $0.0$};
    \foreach \y in {0.5,1.0}
    \draw (0,\y -| origin) node[left] {\small $\y$};
    \foreach \y in {0.5,1.0} {
      \draw[shift={(0,\y -| y)}] (5pt,0pt) -- (0pt,0pt);
      \draw[shift={(0,\y -| x)}] (-5pt,0pt) -- (0pt,0pt);
    }
    \foreach \y in {0.5,1.0}
    {
      \foreach \z in {0.1,0.2,0.3,0.4} {
        \draw[shift={(0,\y -| y)},shift={(0,-0.5)},shift={(0,\z)}]
        (3pt,0pt) -- (0pt,0pt);
        \draw[shift={(0,\y -| x)},shift={(0,-0.5)},shift={(0,\z)}]
        (-3pt,0pt) -- (0pt,0pt);
      }
    }
    \draw (0,0 |- origin) node[below,xshift=2pt] {\small $0$};
    \foreach \x in {5,10,15,20}
    \draw (\x,0 |- origin) node[below] {\small $\x$};
    \foreach \x in {5,10,15} {
      \draw[shift={(\x,0 |- x)}] (0pt,5pt) -- (0pt,0pt);
      \draw[shift={(\x,0 |- y)}] (0pt,-5pt) -- (0pt,0pt);
    }
    \foreach \x in {5,10,15,20}
    {
      \foreach \z in {1,2,3,4} {
        \draw[shift={(\x,0 |- x)},shift={(-5,0)},shift={(\z,0)}]
        (0pt,3pt) -- (0pt,0pt);
        \draw[shift={(\x,0 |- y)},shift={(-5,0)},shift={(\z,0)}]
        (0pt,-3pt) -- (0pt,0pt);
      }
    }
    \draw[thick] (origin) to (x) node[midway,below=10pt] {$t$}
    to (y -| x) to (y) to (origin)
    node[midway,left=20pt] (ylabel)
    {$\displaystyle \frac{\sigma_t}{\sigma_0}$};
    \tikzstyle{EQ}=[red,mark=o,mark size=2pt,mark options={thin,solid}]
    \tikzstyle{PBE}=[blue,mark=x,mark size=2.5pt,mark options={thin,solid}]
    \tikzstyle{MR}=[orange,mark=triangle,mark size=2.5pt,
    mark options={thin,solid}]
    \tikzstyle{rho1}=[solid,thin]
    \tikzstyle{rho10}=[dashed,thin]
    \tikzstyle{rho100}=[dotted,thick]
    \draw[rho1,EQ] plot file {policy-stats-20-eq-1.table};
    \draw[rho10,EQ] plot file {policy-stats-20-eq-10.table};
    \draw[rho100,EQ] plot file {policy-stats-20-eq-100.table};
    \draw[rho1,PBE] plot file {policy-stats-20-pbe-1.table};
    \draw[rho10,PBE] plot file {policy-stats-20-pbe-10.table};
    \draw[rho100,PBE] plot file {policy-stats-20-pbe-100.table};
    \draw[rho1,MR] plot file {policy-stats-20-fl-1.table};
    \draw[rho10,MR] plot file {policy-stats-20-fl-10.table};
    \draw[rho100,MR] plot file {policy-stats-20-fl-100.table};
    \matrix [draw,fill=white,column 2/.style={anchor=mid west},
    right=10pt,ampersand replacement=\&] at (x |- ylabel)
    {
      \draw[rho1,MR,mark repeat=2,mark phase=2] plot coordinates 
      {([xshift=-7.5pt] 0,0) (0,0) ([xshift=7.5pt] 0,0) };
      \&
      \node {$\pi^{\text{MR}},\ \rho_0=1$};
      \\
      \draw[rho10,MR,mark repeat=2,mark phase=2] plot coordinates 
      {([xshift=-7.5pt] 0,0) (0,0) ([xshift=7.5pt] 0,0) };
      \&
      \node {$\pi^{\text{MR}},\ \rho_0=10$};
      \\
      \draw[rho100,MR,mark repeat=2,mark phase=2] plot coordinates 
      {([xshift=-7.5pt] 0,0) (0,0) ([xshift=7.5pt] 0,0) };
      \&
      \node {$\pi^{\text{MR}},\ \rho_0=100$};
      \\
      \draw[rho1,PBE,mark repeat=2,mark phase=2] plot coordinates 
      {([xshift=-7.5pt] 0,0) (0,0) ([xshift=7.5pt] 0,0) };
      \&
      \node {$\pi^{*},\ \rho_0=1$};
      \\
      \draw[rho10,PBE,mark repeat=2,mark phase=2] plot coordinates 
      {([xshift=-7.5pt] 0,0) (0,0) ([xshift=7.5pt] 0,0) };
      \&
      \node {$\pi^{*},\ \rho_0=10$};
      \\
      \draw[rho100,PBE,mark repeat=2,mark phase=2] plot coordinates 
      {([xshift=-7.5pt] 0,0) (0,0) ([xshift=7.5pt] 0,0) };
      \&
      \node {$\pi^{*},\ \rho_0=100$};
      \\
      \draw[rho1,EQ,mark repeat=2,mark phase=2] plot coordinates 
      {([xshift=-7.5pt] 0,0) (0,0) ([xshift=7.5pt] 0,0) };
      \&
      \node {$\pi^{\text{EQ}},\ \rho_0=1$};
      \\
      \draw[rho10,EQ,mark repeat=2,mark phase=2] plot coordinates 
      {([xshift=-7.5pt] 0,0) (0,0) ([xshift=7.5pt] 0,0) };
      \&
      \node {$\pi^{\text{EQ}},\ \rho_0=10$};
      \\
      \draw[rho100,EQ,mark repeat=2,mark phase=2] plot coordinates 
      {([xshift=-7.5pt] 0,0) (0,0) ([xshift=7.5pt] 0,0) };
      \&
      \node {$\pi^{\text{EQ}},\ \rho_0=100$};
      \\
    };
  \end{tikzpicture}
}
  \caption{The evolution  of relative uncertainty 
     of the trader's position for the time horizon $T=20$.
     \label{fig:policy-stats}}
\end{figure}

Under any linear policy, the evolution of the relative uncertainty
$\sigma_t/\sigma_0$ over time is deterministic and depends only on the
parameter $\rho_0$. This is because of the fact that
$\sigma_t/\sigma_0=\rho_t/\rho_0$ and the results in
Section~\ref{se:policy-rep}.  In Figure~\ref{fig:policy-stats}, the
evolution of the relative uncertainty of the PBE policy is
illustrated, for different values of $\rho_0$, as compared to the
equipartitioning and minimum revelation policies. In the low relative
volume regime, the relative uncertainty of the PBE policy evolves
similarly to that of the equipartitioning policy. In the high relative
volume regime, very little information is revealed until close to the
end of trading period under the PBE policy. Indeed, the relative
uncertainty between the equilibrium and the minimum revelation
policies are indistinguishable on the scale of
Figure~\ref{fig:policy-stats}, when $\rho_0=10$ or $\rho_0=100$. These
observations are consistent with our results from
Section~\ref{se:rec}.


\subsection{Adaptive Trading}\label{se:dyn-trade}

One important feature of the linear-Gaussian PBE policy is that it is
adaptive in the sense that the trades executed are random quantities
that are dependent on the exogenous, stochastic fluctuations of the
market. This is in contrast to the policies developed in most of the
optimal execution literature. For example, the baseline
equipartitioning policy of \citet{BertsimasLo98} specifies a
deterministic sequence of trades. Static policies have also been
derived under more complicated models
\citep[e.g.,][]{AlmgrenChriss00,HubermanStanzl05,ObizhaevaWang05,
  AlfonsiSchiedSchulz07b}.  However, this behavior is in contrast to
what is observed amongst institutional traders and trading algorithms
that are implemented by practitioners. One justification for adaptive,
price-responsive trading strategies is risk aversion. It has been
observed that optimal policies for certain risk averse objectives
require dynamic trading \citep{Hora06,Almgren06}. Our model provides
another justification: in the presence of asymmetric information and a
strategic adversary, a trader should seek to exploit price fluctuations
so as disguise trading activity.

In order to understand the behavior of linear policies, it is helpful
to decompose them into deterministic and stochastic components.
Suppose that $(\pi,\psi)$ are a pair of linear policies, and that $y_0
= \mu_0 = 0$. Given Definition~\ref{def:linearpolicy} and
Theorem~\ref{th:Gaussian}, it is easy to see that, for each $ 1 \leq t
\leq T$, there exist vectors $\alpha_{\epsilon,t}, \beta_{\epsilon,t},
\gamma_{\epsilon,t} \in \R^t$ and scalars $\alpha_{x_0,t},
\beta_{x_0,t}, \gamma_{x_0,t} \in \R$, each of which depend on
the parameters $\{\sigma_0,\lambda,\sigma_\epsilon\}$ only
through the $\rho_0$, such that
\begin{equation}\label{eq:dynamictrading}
x_t = \alpha_{x_0,t} x_0 + \frac{1}{\lambda}\alpha_{\epsilon,t}^\top \epsilon^t, 
\quad y_t = \beta_{x_0,t} x_0 + \frac{1}{\lambda}\beta_{\epsilon,t}^\top \epsilon^t, 
\quad \mu_t = \gamma_{x_0,t} x_0 + \frac{1}{\lambda}\gamma_{\epsilon,t}^\top \epsilon^t.
\end{equation}
Here, $\epsilon^t = (\epsilon_1,\ldots,\epsilon_t)$ is the vector of
exogenous disturbances up to time $t$. The first terms in
\eqref{eq:dynamictrading} represent deterministic components of the
policy and the second terms represent zero-mean stochastic components
that depend on market price fluctuations.  For the equipartitioning
and minimum revelation policies, the stochastic components are
zero. On the other hand, the Gaussian PBE policy does have non-zero
stochastic components.

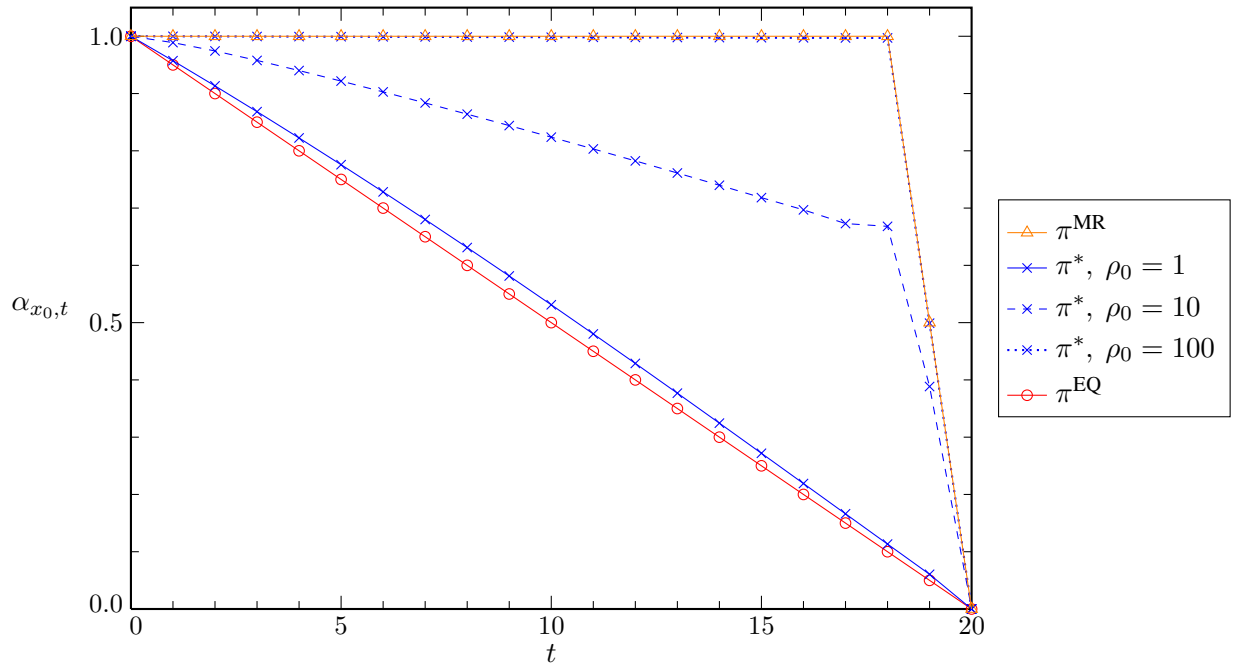
\begin{figure}[htb]
  \centering
  \hidefastcompile{
  \begin{tikzpicture}[x=0.22in,y=3.0in]
    \node[coordinate] at (0,0) (origin) {};
    \node[coordinate] at (20,0 |- origin) (x) {};
    \node[coordinate] at (0,1.05 -| origin) (y) {};
    \draw (0,0 -| origin) node[left,yshift=2pt] {\small $0.0$};
    \foreach \y in {0.5,1.0}
    \draw (0,\y -| origin) node[left] {\small $\y$};
    \foreach \y in {0.5,1.0} {
      \draw[shift={(0,\y -| y)}] (5pt,0pt) -- (0pt,0pt);
      \draw[shift={(0,\y -| x)}] (-5pt,0pt) -- (0pt,0pt);
    }
    \foreach \y in {0.5,1.0}
    {
      \foreach \z in {0.1,0.2,0.3,0.4} {
        \draw[shift={(0,\y -| y)},shift={(0,-0.5)},shift={(0,\z)}]
        (3pt,0pt) -- (0pt,0pt);
        \draw[shift={(0,\y -| x)},shift={(0,-0.5)},shift={(0,\z)}]
        (-3pt,0pt) -- (0pt,0pt);
      }
    }
    \draw (0,0 |- origin) node[below,xshift=2pt] {\small $0$};
    \foreach \x in {5,10,15,20}
    \draw (\x,0 |- origin) node[below] {\small $\x$};
    \foreach \x in {5,10,15} {
      \draw[shift={(\x,0 |- x)}] (0pt,5pt) -- (0pt,0pt);
      \draw[shift={(\x,0 |- y)}] (0pt,-5pt) -- (0pt,0pt);
    }
    \foreach \x in {5,10,15,20}
    {
      \foreach \z in {1,2,3,4} {
        \draw[shift={(\x,0 |- x)},shift={(-5,0)},shift={(\z,0)}]
        (0pt,3pt) -- (0pt,0pt);
        \draw[shift={(\x,0 |- y)},shift={(-5,0)},shift={(\z,0)}]
        (0pt,-3pt) -- (0pt,0pt);
      }
    }
    \draw[thick] (origin) to (x) node[midway,below=10pt] {$t$}
    to (y -| x) to (y) to (origin)
    node[midway,left=20pt] (ylabel)
    {$\alpha_{x_0,t}$};
    \tikzstyle{EQ}=[red,mark=o,mark size=2pt,mark options={thin,solid}]
    \tikzstyle{PBE}=[blue,mark=x,mark size=2.5pt,mark options={thin,solid}]
    \tikzstyle{MR}=[orange,mark=triangle,mark size=2.5pt,
    mark options={thin,solid}]
    \tikzstyle{rho1}=[solid,thin]
    \tikzstyle{rho10}=[dashed,thin]
    \tikzstyle{rho100}=[dotted,thick]
    \draw[rho1,EQ] plot file {policy-deterministic-20-eq-1.table};
    \draw[rho1,PBE] plot file {policy-deterministic-20-pbe-1.table};
    \draw[rho10,PBE] plot file {policy-deterministic-20-pbe-10.table};
    \draw[rho100,PBE] plot file {policy-deterministic-20-pbe-100.table};
    \draw[rho1,MR] plot file {policy-deterministic-20-fl-1.table};
    \matrix [draw,fill=white,column 2/.style={anchor=mid west},
    right=10pt,ampersand replacement=\&] at (x |- ylabel)
    {
      \draw[rho1,MR,mark repeat=2,mark phase=2] plot coordinates 
      {([xshift=-7.5pt] 0,0) (0,0) ([xshift=7.5pt] 0,0) };
      \&
      \node {$\pi^{\text{MR}}$};
      \\
      \draw[rho1,PBE,mark repeat=2,mark phase=2] plot coordinates 
      {([xshift=-7.5pt] 0,0) (0,0) ([xshift=7.5pt] 0,0) };
      \&
      \node {$\pi^{*},\ \rho_0=1$};
      \\
      \draw[rho10,PBE,mark repeat=2,mark phase=2] plot coordinates 
      {([xshift=-7.5pt] 0,0) (0,0) ([xshift=7.5pt] 0,0) };
      \&
      \node {$\pi^{*},\ \rho_0=10$};
      \\
      \draw[rho100,PBE,mark repeat=2,mark phase=2] plot coordinates 
      {([xshift=-7.5pt] 0,0) (0,0) ([xshift=7.5pt] 0,0) };
      \&
      \node {$\pi^{*},\ \rho_0=100$};
      \\
      \draw[rho1,EQ,mark repeat=2,mark phase=2] plot coordinates 
      {([xshift=-7.5pt] 0,0) (0,0) ([xshift=7.5pt] 0,0) };
      \&
      \node {$\pi^{\text{EQ}}$};
      \\
    };
  \end{tikzpicture}
}
  \caption{The deterministic components of trading strategies for the time horizon $T=20$.
    \label{fig:expected-position}}
\end{figure}

Figure~\ref{fig:expected-position} shows the deterministic component
 of the linear-Gaussian PBE versus those of the equipartitioning and minimum
revelation policies. As $\rho_0\tends 0$, the trader
ignores the presence of the arbitrageur and the PBE policy approaches
the equipartitioning policy.  At the other extreme, as
$\rho_0\tends\infty$, in equilibrium the trader seeks to conceal his
activity as much as possible, and hence the PBE policy approaches the
minimum revelation policy.

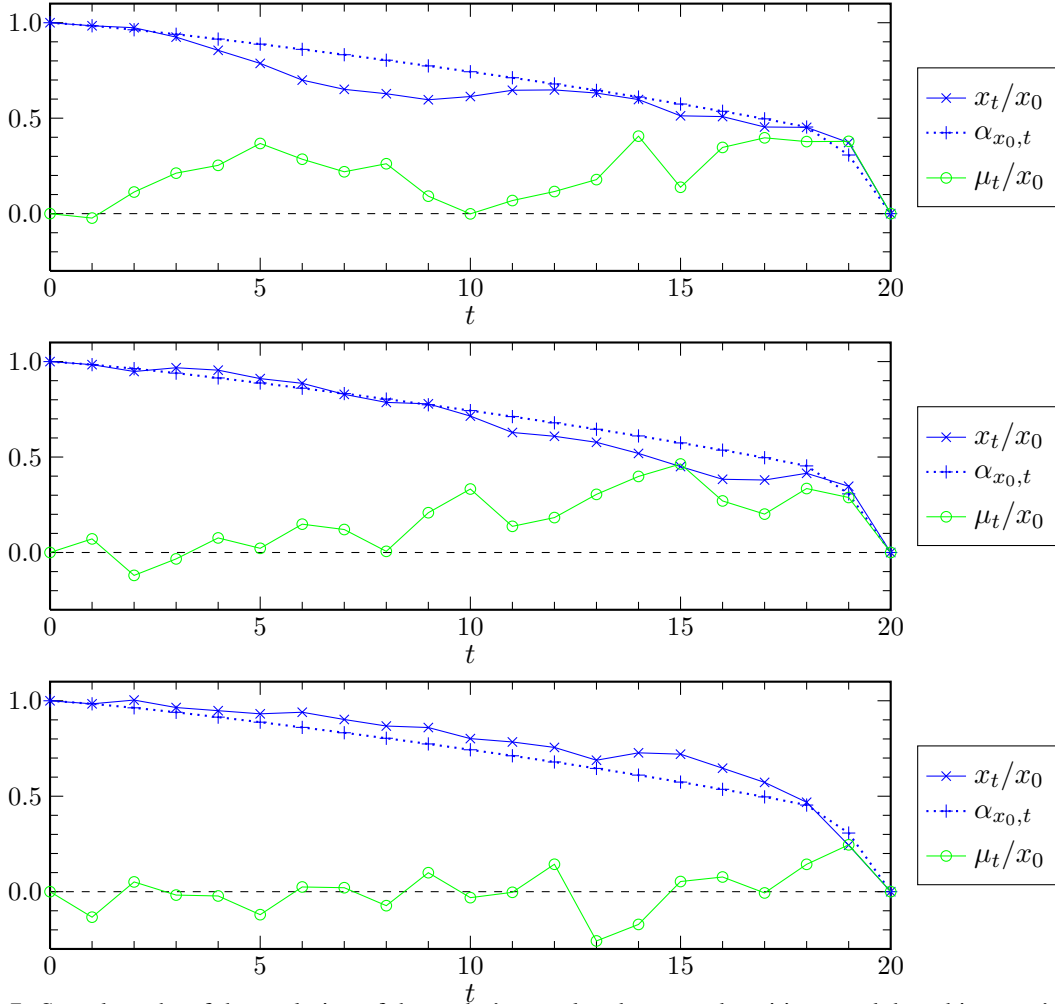
\begin{figure}[htb]
  \centering
  \hidefastcompile{
  \newcommand{\samplepathpic}[1]
  {
    \begin{tikzpicture}[x=0.22in,y=1.0in]
      \node[coordinate] at (0,-0.3) (origin) {};
      \node[coordinate] at (20,0 |- origin) (x) {};
      \node[coordinate] at (0,0.5 -| origin) (midy) {};
      \node[coordinate] at (0,1.1 -| origin) (y) {};
      \foreach \y in {0.0,0.5,1.0}
      \draw (0,\y -| origin) node[left] {\small $\y$};
      \foreach \y in {0.0, 0.5,1.0} {
        \draw[shift={(0,\y -| y)}] (5pt,0pt) -- (0pt,0pt);
        \draw[shift={(0,\y -| x)}] (-5pt,0pt) -- (0pt,0pt);
      }
      \foreach \y in {0.5,1.0}
      {
        \foreach \z in {0.1,0.2,0.3,0.4} {
          \draw[shift={(0,\y -| y)},shift={(0,-0.5)},shift={(0,\z)}]
          (3pt,0pt) -- (0pt,0pt);
          \draw[shift={(0,\y -| x)},shift={(0,-0.5)},shift={(0,\z)}]
          (-3pt,0pt) -- (0pt,0pt);
        }
      }
      \foreach \y in {0.0}
      {
        \foreach \z in {0.3,0.4} {
          \draw[shift={(0,\y -| y)},shift={(0,-0.5)},shift={(0,\z)}]
          (3pt,0pt) -- (0pt,0pt);
          \draw[shift={(0,\y -| x)},shift={(0,-0.5)},shift={(0,\z)}]
          (-3pt,0pt) -- (0pt,0pt);
        }
      }
      \foreach \x in {0,5,10,15,20}
      \draw (\x,0 |- origin) node[below] {\small $\x$};
      \foreach \x in {5,10,15} {
        \draw[shift={(\x,0 |- x)}] (0pt,5pt) -- (0pt,0pt);
        \draw[shift={(\x,0 |- y)}] (0pt,-5pt) -- (0pt,0pt);
      }
      \foreach \x in {5,10,15,20}
      {
        \foreach \z in {1,2,3,4} {
          \draw[shift={(\x,0 |- x)},shift={(-5,0)},shift={(\z,0)}]
          (0pt,3pt) -- (0pt,0pt);
          \draw[shift={(\x,0 |- y)},shift={(-5,0)},shift={(\z,0)}]
          (0pt,-3pt) -- (0pt,0pt);
        }
      }
      \draw[thick] (origin) to (x) node[midway,below=10pt] {$t$}
      to (y -| x) to (y) to (origin)
      node[midway,left=20pt] (ylabel)
      {};
      \draw[dashed] (0,0 -| origin) -- (0,0 -| x);
      \tikzstyle{position}=[blue,thin,solid,mark=x,mark size=2.5pt,
      mark options={thin,solid}]
      \tikzstyle{alpha}=[blue,thick,dotted,mark=+,mark size=2.5pt,
      mark options={thin,solid}]
      \tikzstyle{mu}=[green,thin,solid,mark=o,mark size=2pt,
      mark options={thin,solid}]
      \draw[position] plot file {position-sample-path-20-tr-#1.table};
      \draw[alpha] plot file {position-sample-path-20-alpha-#1.table};
      \draw[mu] plot file {position-sample-path-20-mu-#1.table};
      \matrix [draw,fill=white,column 2/.style={anchor=mid west},
      right=10pt,ampersand replacement=\&] at (x |- ylabel)
      {
        \draw[position,mark repeat=2,mark phase=2] plot coordinates 
        {([xshift=-7.5pt] 0,0) (0,0) ([xshift=7.5pt] 0,0) };
        \pgfmatrixnextcell
        \node {$x_t/x_0$};
        \pgfmatrixendrow
        \draw[alpha,mark repeat=2,mark phase=2] plot coordinates 
        {([xshift=-7.5pt] 0,0) (0,0) ([xshift=7.5pt] 0,0) };
        \pgfmatrixnextcell
        \node {$\alpha_{x_0,t}$};
        \pgfmatrixendrow
        \draw[mu,mark repeat=2,mark phase=2] plot coordinates 
        {([xshift=-7.5pt] 0,0) (0,0) ([xshift=7.5pt] 0,0) };
        \pgfmatrixnextcell
        \node {$\mu_t/x_0$};
        \pgfmatrixendrow
      };
    \end{tikzpicture}
  }
  \samplepathpic{1}
  \smallskip

  \samplepathpic{2}
  \smallskip

  \samplepathpic{3}
}
  \caption{Sample paths of the evolution of the trader's actual and expected positions, and the arbitrageur's mean belief, when $T = 20$, $x_0 = \sigma_0 = 10^5$, $\mu_0 = y_0 = 0$, $\sigma_\epsilon = 0.125$, $\lambda = 10^{-5}$.
    \label{fig:sample-paths}}
\end{figure}

Figure~\ref{fig:sample-paths} illustrates sample paths of the trader's
position under the linear-Gaussian PBE policy. Along each path, the trader
deviates from the deterministic schedule based on the random fluctuations of
the market and how they influence the arbitrageur's beliefs. In general, if
the arbitrageur's estimate of the trader's position becomes more accurate, the
trader accelerates his selling to avoid front-running. On the other hand, if
the arbitrageur is misled as to the trader's position, the trader delays his
selling relative to deterministic schedule.


\section{Extensions}
\label{se:extensions}

In this section, we revisit some of the assumptions in the problem formulation
of Section~\ref{se:model}. At a high level, the main feature of
our model that enables tractability is that, in equilibrium, each agent solves
a linear-quadratic Gaussian control problem. 
This requires that the evolution of the model over time be
described by a linear system 
and that the objectives of the trader and arbitrageur be quadratic functions
that decompose additively over time. As we shall see
shortly, there are a number of extensions of the model one may consider,
incorporating important phenomena such as risk aversion and transient price
impact, that maintain this structure. Such extensions remain tractable and can
be addressed using straightforward adaptations of the techniques we have
developed.

\subsection{Time Horizon}\label{se:trading-horizon}

Our model assumes that the trader begins his liquidation at time $1$
and completes it by time $T$, and that this time interval is common
knowledge. In some instances, public knowledge of the beginning and
end of the liquidation interval might be reasonable since, for
example, this interval will often correspond to a single trading
day. More generally, however, it may be desirable to impose
uncertainty on the part of arbitrageur as to the beginning and end of
the liquidation. Unfortunately, it is not clear how to allow for this
in a tractable fashion in our current framework.

The model further assumes that the arbitrageur must liquidate his position by
time $T+1$. Then, the value function of the arbitrageur at time $T$ with
position $y_T$, is given by $V^*_T(y_T) = - \lambda y_T^2$. This was used in
\eqref{eq:terminal-V}--\eqref{eq:terminal-U} to determine the value functions
$U^*_{T-1}$ and $V^*_{T-1}$, which form the base case of the backward
induction.  This assumption can easily be relaxed. For example, suppose that
the arbitrageur has $T_a$ additional trading periods. It is easy to see that,
after time $T$, the arbitrageur will optimally equipartition over the
remaining $T_a$ periods. Therefore the value of a position $y_T$ at time $T$
will take the form $V^*_{T}(y_T) = -\lambda \frac{T_a + 1}{2T_a} y_T^2$,
following the analysis in \citet{BertsimasLo98}. So long as $V^*_T$ is a
quadratic function, our discussion in Sections~\ref{se:dp-analysis} and
\ref{se:algorithm} carries through, with a different choice of terminal value
functions. 


\subsection{Risk Aversion}\label{se:risk-aversion}

Our model assumes that both the trader and arbitrageur are
risk-neutral. One way to account for risk aversion is to follow the 
approach suggested by \citet{Hora06}. In particular, we could assume that, for example,
the trader seeks to optimize the objective function
\[
\E\left[\left.
    \sum_{\tau=0}^{T-1} \left\{
      \Delta p_{\tau+1} x_\tau
      - 
      \frac{\eta}{2}\big(\Delta p_{\tau+1} x_\tau
        - \E[\Delta p_{\tau+1} x_\tau\ |\ x_\tau,y_\tau,\phi_\tau]\big)^2
      - \zeta x_\tau^2
      \right\}
\ \right|\ x_0,y_0,\phi_0 \right],
\]
The second term in the sum penalizes for variance in revenue in each 
time period, with $\eta \geq 0$ capturing the degree of risk aversion.  This final term
represents a per stage holding cost, with the parameter $\zeta\geq 0$
expressing the degree to which the trader would prefer to execute
sooner rather than later. The risk neutral case previously considered
corresponds to the choice of $\eta = \zeta = 0$. For any nonnegative 
parameter choices, the objective remains a time
separable positive definite quadratic function. Hence, the methods of
Sections~\ref{se:dp-analysis} and \ref{se:algorithm} can be suitably
adapted.

\subsection{Price Impact \& Price Dynamics}\label{se:price-impact}

Our model assumed permanent and linear price impact. Empirically, it
has been observed that transient price impact is a significant
component of price dynamics, and it is important to account for this
in the design of execution strategies.

More generally, our analysis applies when there is some collection of
state variables (for example, $\{x_t, y_t, \mu_t\}$) that evolve as a
linear dynamical system with Gaussian disturbances, and where changes
in price are linear in the state variables.  In order to
incorporate transient price impact, assume that prices evolve
according to
\begin{equation}\label{eq:temp-impact}
\begin{split}
p_{t} & = p_{0} + \underbrace{\lambda \sum_{\tau=1}^{t}
(u_\tau + v_\tau + z_\tau)}_{\text{permanent price impact}} 
+ \underbrace{\gamma \sum_{\tau=1}^{t} \alpha^{t-\tau} (u_\tau + v_\tau + z_\tau) 
}_{\text{transient price impact}}.
\end{split}
\end{equation}
Here, $u_\tau$ and $v_\tau$ are the trades of the trader and arbitrageur,
respectively, as time $\tau$. In place of the exogenous noise term in the
original price dynamics \eqref{eq:p-evolve}, $z_\tau$ is an IID
$N(0,\sigma_z^2)$ random variable representing the quantity of noise trades at
time $\tau$. The second term in \eqref{eq:temp-impact} captures a permanent,
linear price impact with sensitivity $\lambda \geq 0$. The final term
represents a transient, linear price impact with sensitivity $\gamma \geq 0$
and recovery rate $\alpha \in [0,1)$.

These price dynamics can be rewritten as 
\[
p_t = p_{t-1} + (\lambda + \gamma) (u_t + v_t + z_t) 
- \gamma (1 - \alpha) s_{t-1},
\]
where $s_{t}$ is defined to be geometrically weighted total order flow
\[
s_{t} \defeq \sum_{\tau=0}^t \alpha^{t-\tau} (u_\tau + v_\tau + z_\tau) 
= \alpha s_{t-1} + (u_t + v_t + z_t).
\]
Now, suppose that the trader's decision $u_t$ is a linear function of $\{
x_{t-1}, y_{t-1}, \mu_{t-1}, s_{t-1}\}$, and the arbitrageur's decision $v_t$
is a linear function of $\{ y_{t-1}, \mu_{t-1}, s_{t-1}\}$. Then, it will be
the case that $\{ x_t, y_t, \mu_t, s_t \}$ evolve as a linear dynamical
system, and that the price changes are linear in these state
variables. Therefore, the analysis in Sections~\ref{se:dp-analysis} and
\ref{se:algorithm} can be suitably modified and repeated, with an augmented
state space. Note that, since $s_t$ is a function of only of the {\em total}
quantities traded at times up to $t$, it is reasonable to assume that this is
public knowledge known to both the trader and arbitrageur.

Other aspects of more complicated price dynamics can also be incorporated via
such state augmentation. For example, one may consider linear factor models or
other otherwise add exogenous explanatory variables to the evolution of
prices, so long as the dependencies are linear. Similarly, models that
incorporate drift in the price process, such as short term momentum or mean
reversion, can be considered.

\subsection{Parameterized Policies}

Beyond solving specific classes of models, results from the optimal execution
literature offer useful guidance on how to structure parameterized execution
policies that can be effective even if modeling assumptions are not entirely
valid.  In this vein, concepts we have developed can enhance parameterized
policies that one might design based on prior literature.

For example, consider designing an execution system which begins the trading
day with a position that must be liquidated by the end of that trading day.  A
number of models previously considered in the literature result in
deterministic linear policies \citep[see,
e.g.,][]{BertsimasLo98,ObizhaevaWang05,AlfonsiSchiedSchulz07a}.  In
particular, for each $t$th time period during the course of the day, there is
a parameter $a_t$ that indicates the fraction of the position to sell during
that time period.  These parameters $a_0, \ldots, a_{T-1}$ depend on
asset-specific characteristics such as volatility and market impact model
parameters.

Modeling assumptions often do not match reality. As such, it is useful to add
flexibility by parametrizing the execution policy. For example, we might
employ a policy that sells a fraction $\theta_t a_t$ of the position during
each $t$th time period, where $\theta_0,\ldots,\theta_{T-1}$ are
asset-independent parameters. Then, these parameters can be tuned based on
experience from trading all assets. It is important that the number of
parameters does not scale with the number of assets, because we would then be
unlikely to have a sufficient amount of data to tune parameters. In this
regard, the way $a_0, \ldots, a_{T-1}$ capture variations across assets is
critical to the design of an effective parametrization.

Our work motivates a generalized class of parameterized policies that adapt
trades as price movements are observed. Our model is optimized by an execution
strategy with three sequences of coefficients: $\{a_{x,t}, a_{y,t},
a_{\mu,t}\ |\ t=0,\ldots,T-1\}$. By simulating arbitrageur activity over the
course of the day and applying these coefficients appropriately, we produce a
sequence of trades that adapt to price fluctuations. Similarly with the case
of a deterministic policy, we can introduce parameters $\{\theta_{x,t},
\theta_{y,t}, \theta_{\mu,t}\ |\ t=0,\ldots,T-1\}$ that scale the policy
coefficients, and tune these parameters based on experience. Once again, these
parameters are asset-independent while the coefficients $\{a_{x,t}, a_{y,t},
a_{\mu,t}\ |\ t=0,\ldots,T-1\}$ capture dependence of the policy on
asset-specific characteristics such as volatility and market impact model
parameters.

\section{Conclusion}
\label{se:conclusion}

Our model captures strategic interactions between a trader aiming to liquidate
a position and an arbitrageur trying to detect and profit from the trader's
activity. The algorithm we have developed computes Gaussian perfect Bayesian
equilibrium behavior. It is interesting that the resulting trader policy takes
on such a simple form: the number of shares to liquidate at time $t$ is linear
in the trader's position $x_{t-1}$, the arbitrageur's position $y_{t-1}$ and
the arbitrageur's estimate $\mu_{t-1}$ of $x_{t-1}$. The coefficients of the
policy depend only on the relative volume parameter $\rho_0$, which quantifies
the magnitude of the trader's position relative to the typical market
activity, and the time horizon $T$. This policy offers useful guidance beyond
what has been derived in models that do not account for arbitrageur behavior.
In the absence of an arbitrageur, it is optimal to trade equal amounts over
each time period, which corresponds to a policy that is linear in $x_{t-1}$.
The difference in the PBE policy stems from its accounting of the
arbitrageur's inference process. In particular, the policy reduces information
revealed to the arbitrageur by delaying trades and takes advantage of
situations where the arbitrageur has been misled by unusual market activity.

Our model represents a starting point for the study of game theoretic
behavior in trade execution. It has an admittedly simple structure,
and this allows for a tractable analysis that highlights the
importance of information signaling. There are a number of extensions
to this model that are possible, however, and that warrant further
discussion:
\begin{enumerate}

\item {\bf (Flexible Time Horizon)} We assume finite time horizons
  $T$ and $T+1$ for the trader and arbitrageur, respectively. 
  The choice of time horizon has
  an impact on the resulting equilibrium policies, and there are
  clearly end-of-horizon effects in the policies computed in
  Section~\ref{se:computation}. To some extent it seems artificial to
  impose a fixed time horizon as an exogenous restriction on
  behavior. Fixed horizon models preclude the trader from delaying
  liquidation beyond the horizon even if this can yield significant
  benefits, for example. A better model would be to consider an
  infinite horizon game, where risk aversion provides the motivation
  for liquidating a position sooner rather than later.

\item {\bf (Uncertain Trader)} In our model, we assume that the
  arbitrageur is uncertain of the trader's position, but that the
  trader knows everything. A more realistic model would allow for
  uncertainty on the part of the trader as well, and would allow for
  the arbitrageur to mislead the trader.

\item {\bf (Multi-player Games)} Our model restricts to a single trader
  and arbitrageur. A natural extension would be to consider multiple
  traders and arbitrageurs that are uncertain about each others'
  positions and must compete in the marketplace as they unwind. Such a
  generalized model could be useful for analysis of important
  liquidity issues such as those arising from the credit crunch of
  2007.
\end{enumerate}

Also of interest are the potential empirical implications of the model. If we
make the assumption that the trade execution horizon is a single day, the
observations in Section~\ref{se:computation} suggest particular patterns for
intraday volume. For example, if $\rho_0$ is large, the volume traded should
be much higher near the end of the day then at other times. Similarly, the
structure of the equilibrium trading policies for the trader and arbitrageur
will generate specific, time-varying auto-correlation in the increments of the
price process. Formulating tests of such empirical predictions in any
interesting area for future research.

Finally, beyond the immediate context of our model, there are many
directions worth exploring.  One important avenue is to factor data
beyond price into the execution strategy.  For example, volume data
may play a significant role in the arbitrageur's inference, in which
case it should also influence execution decisions.  Limit order book
data may also be relevant.  Developing tractable models that account
for such data remains a challenge.  One initiative to incorporate
limit order book data into the decision process is presented by
\citet{Nevmyvaka06}.

\section*{Acknowledgments}

The authors are grateful to the anonymous reviewers for helpful comments and
corrections. The first author wishes to thank Mark Broadie and Gabriel
Weintraub for useful discussions. The second author wishes to thank the
Samsung Scholarship Foundation for financial support. This research was
conducted while the third author was visiting the Faculty of Commerce and
Accountancy at Chulalongkorn University and supported by the Chin Sophonpanich
Foundation Fund.

{\small
\singlespacing
\bibliography{execution}
}

\appendix

\section{Proofs}

\begin{th:Gaussian}
  If the belief distribution $\phi_{t-1}$ at time is Gaussian, and the
  arbitrageur assumes that the trader's policy $\hat{\pi}_t$ is linear with
  $\hat{\pi}_t(x_{t-1},y_{t-1},\phi_{t-1}) = \hat{a}_{x,t}^{\rho_{t-1}}
  x_{t-1} + \hat{a}_{y,t}^{\rho_{t-1}} y_{t-1} + \hat{a}_{\mu,t}^{\rho_{t-1}}
  \mu_{t-1}$, then the belief distribution $\phi_t$ is also Gaussian. The mean
  $\mu_t$ is a linear function of $y_{t-1}$, $\mu_{t-1}$, and the observed
  price change $\Delta p_t$, with coefficients that are deterministic
  functions of the scaled variance $\rho_{t-1}$. The scaled variance $\rho_t$
  evolves according to
\[
\rho_t^2 = \left(1+\hat{a}^{\rho_{t-1}}_{x,t}\right)^2
\left(\frac{1}{\rho_{t-1}^2}
+ (\hat{a}^{\rho_{t-1}}_{x,t})^2
\right)^{-1}.
\]
In particular, $\rho_t$ is a deterministic function of $\rho_{t-1}$.
\end{th:Gaussian}
\begin{proof}
Set $\{K_{t-1},h_{t-1}\}$ to be the
information form parameters for the Gaussian distribution
$\phi_{t-1}$, so that
\[
K_{t-1} \defeq 1/\sigma^{2}_{t-1}, \qquad \text{and}\qquad 
h_{t-1} \defeq \mu_{t-1}/\sigma^2_{t-1}.
\]
Define $\phi_{t-1}^+$ to be the distribution of $x_{t-1}$ conditioned
on all information seen by the arbitrageur at times up to and
including $t$. That is,
\[
\phi_{t-1}^+(S) \defeq \Pr\big(x_{t-1} \in S \ |\ \phi_{t-1}, y_{t-1}, 
\lambda(\hat{\pi}_t(x_{t-1}, y_{t-1}, \phi_{t-1}) + v_t) 
+ \epsilon_t = \Delta p_t\big),
\]
where $\Delta p_t$ is the price change observed at time $t$.  By
Bayes' rule, this distribution has density
\[
\begin{split}
\phi^+_{t-1}(dx) & \propto
\phi_{t-1}(dx) \,
\exp\left( - \frac{\bigl(\Delta p_{t} 
- \lambda (\pi_{t}(x,y_{t-1},\phi_{t-1})
+ \psi_{t}(y_{t-1},\phi_{t-1}))
\bigr)^2}{2 \sigma^2_\epsilon}
\right)
\\
& \propto
\exp\Biggl(
-\tfrac{1}{2} 
\left(K_{t-1} 
+ \frac{\lambda^2 (\hat{a}_{x,t}^{\rho_{t-1}})^2}{\sigma^2_\epsilon}
\right) 
x^2
\\
& \quad\quad\quad\quad
+
\left(
h_{t-1} 
+
\frac{\lambda \bigl(\Delta p_{t}- \lambda (\hat{a}_{y,t}^{\rho_{t-1}} y_{t-1} 
+ \hat{a}_{\mu,t}^{\rho_{t-1}}  \mu_{t-1}
+ \psi_{t})\bigr)
\hat{a}_{x,t}}{\sigma^2_\epsilon}
\right)
x 
\Biggr)\, dx.
\end{split}
\]
Thus, $\phi^+_{t-1}$ is a Gaussian distribution, with variance
\[
\left(K_{t-1} 
+ \frac{\lambda^2 (\hat{a}_{x,t}^{\rho_{t-1}})^2}{\sigma^2_\epsilon}
\right)^{-1},
\]
and mean
\[
\left(K_{t-1} 
+ \frac{\lambda^2 (\hat{a}_{x,t}^{\rho_{t-1}})^2}{\sigma^2_\epsilon}
\right)^{-1}
\left(
h_{t-1} 
+
\frac{\lambda \bigl(\Delta p_{t}- \lambda (\hat{a}_{y,t}^{\rho_{t-1}} y_{t-1} 
+ \hat{a}_{\mu,t}^{\rho_{t-1}}  \mu_{t-1}
+ \psi_{t})\bigr)
\hat{a}_{x,t}}{\sigma^2_\epsilon}
\right).
\]

Now, note that
\[
x_{t} = x_{t-1} + \hat{\pi}_{t}(x_{t-1},y_{t-1},\phi_{t-1})
= (1 + \hat{a}_{x,t}^{\rho_{t-1}}) x_{t-1} + \hat{a}_{y,t}^{\rho_{t-1}} y_{t-1} + \hat{a}_{\mu,t}^{\rho_{t-1}} \mu_{t-1}.
\]
Then, $\phi_t$ is also a Gaussian distribution, with variance
\begin{equation}\label{eq:var-update}
\sigma^2_t  =
(1 + \hat{a}_{x,t}^{\rho_{t-1}})^2
\left(K_{t-1} 
+ \frac{\lambda^2 (\hat{a}_{x,t}^{\rho_{t-1}})^2}{\sigma^2_\epsilon}
\right)^{-1}
=
(1 + \hat{a}_{x,t})^2
\left(\frac{1}{\sigma_{t-1}^2}
+ \frac{\lambda^2 (\hat{a}_{x,t}^{\rho_{t-1}})^2}{\sigma^2_\epsilon}
\right)^{-1},
\end{equation}
and mean
\begin{equation}\label{eq:mean-update}
\begin{split}
\mu_t &=
\hat{a}_{y,t}^{\rho_{t-1}} y_{t-1} + \hat{a}_{\mu,t}^{\rho_{t-1}} \mu_{t-1}
\\
& \quad
+
(1 + \hat{a}_{x,t}^{\rho_{t-1}})
\frac{
\mu_{t-1}/\rho^2_{t-1}
+
\bigl(\Delta p_{t}/\lambda - \hat{a}_{y,t}^{\rho_{t-1}} y_{t-1}
- \hat{a}_{\mu,t}^{\rho_{t-1}} \mu_{t-1}
- \psi_{t}\bigr)
\hat{a}_{x,t}
}
{1/\rho^2_{t-1}
+ (\hat{a}_{x,t}^{\rho_{t-1}})^2}.
\end{split}
\end{equation}
The conclusions of the theorem immediately follow.
\end{proof}

In order to prove Theorems~\ref{th:tqd-aqd}--\ref{th:linear-pi-psi}, it is
necessary to explicitly evaluate the operator $F^{(\psi_t, \pi_t)}_{u_t}$
applied to quadratic functions of $\{x_t,y_t,\mu_t\}$ and the operator
$G_{v_t}^{\pi_t}$ applied to quadratic functions of $\{y_t,\mu_t\}$. The
following lemma is helpful for this purpose, as it provides expressions for
the expectation of $\mu_t$ and $\mu_t^2$ under various distributions.

\begin{lemma}\label{le:expmu}
  Assume that the the policies $\psi_t$ and
  $\pi_t$ are  linear with
\begin{gather*}
\pi_t(x_{t-1},y_{t-1},\phi_{t-1}) 
 = a_{x,t}^{\rho_{t-1}} x_{t-1} 
+ a_{y,t}^{\rho_{t-1}} y_{t-1} + a_{\mu,t}^{\rho_{t-1}} \mu_{t-1},
\\
\psi_t(y_{t-1},\phi_{t-1}) 
= b^{\rho_{t-1}}_{y,t} y_{t-1} + b^{\rho_{t-1}}_{\mu,t} \mu_{t-1}.
\end{gather*}
Define
\[
\gamma^{\rho_{t-1}}_t \defeq 
\frac{1 + a^{\rho_{t-1}}_{x,t}}{1/\rho^2_{t-1} 
+ (a^{\rho_{t-1}}_{x,t})^2}.
\]
Then,
\begin{subequations}
\begin{gather}
\begin{split}
\E^{(\psi_t, \pi_t)}_{u_t}
\left[\left. \mu_t\ \right|\ x_{t-1}, y_{t-1}, \phi_{t-1}\right] 
&
= 
a^{\rho_{t-1}}_{y,t} y_{t-1} + a^{\rho_{t-1}}_{\mu,t} \mu_{t-1}
+ \gamma^{\rho_{t-1}}_t \mu_{t-1} / \rho^2_{t-1}
\\
& 
\quad
+ \gamma^{\rho_{t-1}}_t
a^{\rho_{t-1}}_{x,t}
\left(u_t - a^{\rho_{t-1}}_{y,t} y_{t-1} 
- a^{\rho_{t-1}}_{\mu,t} \mu_{t-1}\right),
\end{split}\label{eq:expmu-trader1}
\\[\baselineskip]
\Var^{(\psi_t, \pi_t)}_{u_t}
\left[\left.\mu_t\ \right|\ x_{t-1}, y_{t-1}, \phi_{t-1}\right] 
=
\left(\gamma^{\rho_{t-1}}_t a^{\rho_{t-1}}_{x,t} 
\sigma_\epsilon / \lambda\right)^2,
\\[\baselineskip]
\begin{split}
\E^{(\psi_t, \pi_t)}_{u_t}
\left[\left.\mu^2_t\ \right|\ x_{t-1}, y_{t-1}, \phi_{t-1}\right] 
& 
= \Var^{(\psi_t, \pi_t)}_{u_t}\left[\left.\mu_t\ \right|\ 
x_{t-1}, y_{t-1}, \phi_{t-1}\right]
\\
& \quad
+ \left(\E^{(\psi_t, \pi_t)}_{u_t}
\left[\left.\mu_t\ \right|\ x_{t-1}, y_{t-1},\phi_{t-1}\right]\right)^2,
\end{split}\label{eq:expmu-trader2}
\\[\baselineskip]
\E^{\pi_t}_{v_t}\left[\left.\mu_t\ \right|\ y_{t-1},\phi_{t-1}\right] 
= a^{\rho_{t-1}}_{y,t} y_{t-1} + 
(1+a^{\rho_{t-1}}_{x,t}+a^{\rho_{t-1}}_{\mu,t}) \mu_{t-1}, 
\label{eq:expmu-arb1}
\\[\baselineskip]
\Var^{\pi_t}_{v_t}\left[\left.\mu_t\ \right|\ y_{t-1},\phi_{t-1}\right] 
= \left(
\gamma^{\rho_{t-1}}_t
a^{\rho_{t-1}}_{x,t}
\sigma_\epsilon / \lambda
\right)^2 
\left(1+\left(a^{\rho_{t-1}}_{x,t}\right)^2 \rho^2_{t-1}\right),
\\[\baselineskip]
\E^{\pi_t}_{v_t}\left[\left.\mu^2_t\ \right|\ y_{t-1},\phi_{t-1}\right] 
= \Var^{\pi_t}_{v_t}\left[\left.\mu_t\ \right|\ y_{t-1},\phi_{t-1}\right] 
+ \left(\E^{\pi_t}_{v_t}
\left[\left.\mu_t\ \right|\ y_{t-1},\phi_{t-1}\right]\right)^2.
\label{eq:expmu-arb2}
\end{gather}
\end{subequations}
\end{lemma}
\begin{proof}
  The lemma follows directly from taking expectations of the mean
  update equation \eqref{eq:mean-update}.
\end{proof}

\begin{th:tqd-aqd}
  If $U^*_t$ is TQD and $V^*_t$ is AQD, and
  Step~\ref{algline:lgpbe:pipsi} of Algorithm~\ref{alg:lgpbe} produces
  a linear pair $(\pi_t^*,\psi_t^*)$, then $U^*_{t-1}$ and
  $V^*_{t-1}$, defined by Step~\ref{algline:lgpbe:UV} of
  Algorithm~\ref{alg:lgpbe} are TQD and AQD, respectively.
\end{th:tqd-aqd}
\begin{proof}
Suppose that
\[
\begin{split}
V^*_t(y_t,\phi_t) 
&= - \lambda\left(\tfrac{1}{2} d^{\rho_t}_{yy,t} y_t^2 
+ \tfrac{1}{2} d^{\rho_t}_{\mu\mu,t} \mu_t^2
+ d^{\rho_t}_{y\mu,t} y_t \mu_t
-  \frac{\sigma_\epsilon^2}{\lambda^2} d^{\rho_t}_{0,t}\right),
\\
\pi^*_t(x_{t-1},y_{t-1},\phi_{t-1})
&
= a^{\rho_{t-1}}_{x,t} x_{t-1} + a^{\rho_{t-1}}_{y,t} y_{t-1}
+ a^{\rho_{t-1}}_{\mu,t} \mu_{t-1}, 
\\
\psi^*_t(y_{t-1},\phi_{t-1}) 
&
= b^{\rho_{t-1}}_{y,t} y_{t-1}+ b^{\rho_{t-1}}_{\mu,t} \mu_{t-1}.
\end{split}
\]
If the trader uses the policy $\pi^*_t$ and the arbitrageur uses the
policy $\psi^*_t$, we have
\[
\begin{split}
u_t & = 
a^{\rho_{t-1}}_{x,t} x_{t-1} + a^{\rho_{t-1}}_{y,t} y_{t-1}
+ a^{\rho_{t-1}}_{\mu,t} \mu_{t-1},  \\
v_t & =
 b^{\rho_{t-1}}_{y,t} y_{t-1}+ b^{\rho_{t-1}}_{\mu,t} \mu_{t-1}, \\
y_t & = y_{t-1} +  b^{\rho_{t-1}}_{y,t} y_{t-1}+ b^{\rho_{t-1}}_{\mu,t} \mu_{t-1}.
\end{split}
\]
Using these facts, Theorem~\ref{th:Gaussian}, and
\eqref{eq:expmu-arb1}--\eqref{eq:expmu-arb2} from Lemma~\ref{le:expmu}, we can
explicitly compute
\[
\begin{split}
V^*_{t-1}(y_{t-1},\phi_{t-1}) & =
\left(G_{\psi_t^*}^{\pi_t^*} V\right)(y_{t-1},\phi_{t-1}) 
\\
& = \E^{\pi_t^*}_{\psi_t^*}\left[\lambda (u_t + v_t) y_{t-1} + V^*_t(y_t,\phi_t) 
\ \Big| \ y_{t-1}, \phi_{t-1}\right] \\
& = - \lambda\left(\tfrac{1}{2} d^{\rho_{t-1}}_{yy,t-1} y_t^2 
+ \tfrac{1}{2} d^{\rho_{t-1}}_{\mu\mu,t-1} \mu_t^2
+ d^{\rho_{t-1}}_{y\mu,t-1} y_t \mu_t
-  \frac{\sigma_\epsilon^2}{\lambda^2} d^{\rho_{t-1}}_{0,t-1}\right),
\end{split}
\]
where
\[
\begin{split}
\rho_t^2 & = \left(1+\hat{a}^{\rho_{t-1}}_{x,t}\right)^2
\left(\frac{1}{\rho_{t-1}^2}
+ (\hat{a}^{\rho_{t-1}}_{x,t})^2
\right)^{-1}, \\
d^{\rho_{t-1}}_{yy,t-1} &
= \left(d^{\rho_{t}}_{yy,t} - \frac{(d^{\rho_{t}}_{y\mu,t})^2}{d^{\rho_{t}}_{yy,t}} \right) (a^{\rho_{t}}_{y,t})^2
+ 2 \left(\frac{d^{\rho_{t}}_{y\mu,t}}{d^{\rho_{t}}_{yy,t}} -1\right)a^{\rho_{t}}_{y,t} - \frac{1}{d^{\rho_{t}}_{yy,t}} + 2,
\\
d^{\rho_{t-1}}_{y\mu,t-1} &
= -a^{\rho_{t}}_{\mu,t} - a^{\rho_{t}}_{x,t} + \left(\frac{d^{\rho_{t}}_{y\mu,t}}{d^{\rho_{t}}_{yy,t}} + \left(d^{\rho_{t}}_{\mu\mu,t} -\frac{(d^{\rho_{t}}_{y\mu,t})^2}{d^{\rho_{t}}_{yy,t}}\right)a^{\rho_{t}}_{y,t} \right) (1+a^{\rho_{t}}_{x,t}+a^{\rho_{t}}_{y,t}),
\\
d^{\rho_{t-1}}_{\mu\mu,t-1} &
= \left(d^{\rho_{t}}_{\mu\mu,t} -\frac{(d^{\rho_{t}}_{y\mu,t})^2}{d^{\rho_{t}}_{yy,t}} \right) (1+a^{\rho_{t}}_{x,t}+a^{\rho_{t}}_{y,t})^2,\\
d^{\rho_{t-1}}_{0,t-1} &
= d^{\rho_{t}}_{0,t} + \frac{d^{\rho_{t}}_{\mu\mu,t}}{2} \left(a^{\rho_{t}}_{x,t} \gamma^{\rho_{t-1}}_{t} \frac{\sigma_\epsilon}{\lambda}\right)^2 \left(1 + (\rho_{t-1} a^{\rho_{t}}_{x,t})^2 \right).
\end{split}
\]
Therefore, $V^*_{t-1}$ is AQD. Similarly, we can check that $U^*_{t-1}$ is TQD.
\end{proof}

\begin{th:sufficientcondition}
  Suppose that $U^*_t$ and $V^*_t$ and TQD/AQD value functions
  specified by \eqref{eq:TQD}--\eqref{eq:AQD}, and
  $(\pi^*_t,\psi^*_t)$ are linear policies specified by
  \eqref{eq:linear-pi}--\eqref{eq:linear-psi}.  Assume that, for all
  $\rho_{t-1}$, the policy coefficients satisfy the first order
  conditions
  \begin{gather}
    \begin{split}
      0 & =
      \big(\rho_t^2 c^{\rho_t}_{\mu\mu,t} + 2 \rho_t c^{\rho_t}_{x\mu,t} + c^{\rho_t}_{xx,t}\big)
      \big(a^{\rho_{t-1}}_{x,t}\big)^3 + \big(3c^{\rho_t}_{xx,t} + 3\rho_t
      c^{\rho_t}_{x\mu,t} - 1\big)\big(a^{\rho_{t-1}}_{x,t}\big)^2 \\
      & \quad + \big(3c^{\rho_t}_{xx,t} + \rho_t c^{\rho_t}_{x\mu,t} - 2\big)
      a^{\rho_{t-1}}_{x,t} + c^{\rho_t}_{xx,t} - 1, \\
    \end{split}\label{eq:simple-start2}
    \\
    a^{\rho_{t-1}}_{y,t} 
    = 
    -\frac{\big(b^{\rho_{t-1}}_{y,t}+1\big) \big(c^{\rho_t}_{xy,t}+\alpha_t
      c^{\rho_t}_{y\mu,t}\big)}{c^{\rho_t}_{xx,t}+(\alpha_t +1) 
      c^{\rho_t}_{x \mu,t}+\alpha_t c^{\rho_t}_{\mu \mu,t}},\label{eq:recursion-ay2} 
    \\
    a^{\rho_{t-1}}_{\mu,t} 
    =
    -\frac{a^{\rho_{t-1}}_{x,t} b^{\rho_{t-1}}_{\mu,t}
      \big(c^{\rho_t}_{xy,t}+\alpha_t c^{\rho_t}_{y \mu,t}\big)+\alpha_t 
      \big(c^{\rho_t}_{x \mu,t} + 
      \alpha_t c^{\rho_t}_{\mu\mu,t}\big)/\rho_{t-1}^2}
    {a^{\rho_{t-1}}_{x,t} \big(c^{\rho_t}_{xx,t}+(\alpha_t +1) 
      c^{\rho_t}_{x \mu,t}+\alpha_t c^{\rho_t}_{\mu \mu,t}\big)},\label{eq:recursion-amu2} 
    \\
    b^{\rho_{t-1}}_{y,t} 
    =
    \frac{1-d^{\rho_t}_{y\mu,t}a^{\rho_{t-1}}_{y,t}}{d^{\rho_t}_{yy,t}} - 1,
    \qquad
    b^{\rho_{t-1}}_{\mu,t} 
    =
    -\frac{(1+a^{\rho_{t-1}}_{\mu,t} 
      + a^{\rho_{t-1}}_{x,t})d^{\rho_t}_{y\mu,t}}{d^{\rho_t}_{yy,t}},\label{eq:recursion-by-bmu2}
  \end{gather}
  and the second order conditions
  \begin{equation}\label{eq:secondordercondition2}
    c^{\rho_t}_{xx,t}+(\alpha_t +1) 
    c^{\rho_t}_{x \mu,t}+\alpha_t c^{\rho_t}_{\mu \mu,t} > 0,
    \qquad
    d^{\rho_t}_{yy,t} > 0,
  \end{equation}
  where the quantities $\alpha_t$ and $\rho_t$ satisfy
  \begin{equation}
    \label{eq:simple-end2}
    \alpha_t = \frac{a^{\rho_{t-1}}_{x,t}\big(1 + a^{\rho_{t-1}}_{x,t}\big)}{1/\rho^2_{t-1} + \big(a^{\rho_{t-1}}_{x,t}\big)^2}, 
    \qquad
    \rho_t^2 = \left(1+a^{\rho_{t-1}}_{x,t}\right)^2
    \left(\frac{1}{\rho_{t-1}^2}
      + \big(a^{\rho_{t-1}}_{x,t}\big)^2
    \right)^{-1}.
  \end{equation}
  Then, $(\pi^*_t,\psi^*_t)$ satisfy the single-stage equilibrium conditions
  \begin{align}\label{eq:singlestageproblemU2}
    \pi^*_t(x_{t-1},y_{t-1},\phi_{t-1}) & \in
    \argmax_{u_t} \left(F^{(\psi^*_t, \pi^*_t)}_{u_t}
      U^*_t\right)(x_{t-1},y_{t-1},\phi_{t-1}),
    \\
    \label{eq:singlestageproblemV2}
    \psi^*_t(y_{t-1},\phi_{t-1}) & \in 
    \argmax_{v_t} \left(G_{v_t}^{\pi^*_t}
      V^*_t\right)(y_{t-1},\phi_{t-1}),
  \end{align}
  for all $x_{t-1}$, $y_{t-1}$, and Gaussian $\phi_{t-1}$.
\end{th:sufficientcondition}
\begin{proof}
  As we will discuss in the proof of Theorem \ref{th:linear-pi-psi}, the
  optimizing value $u_t^*$ in \eqref{eq:singlestageproblemU2} is a linear
  function of $x_{t-1}$, $y_{t-1}$ and $z_{t-1}$, whose coefficients depend on
  $\{a_{x,t}^{\rho_{t-1}},a_{y,t}^{\rho_{t-1}},a_{\mu,t}^{\rho_{t-1}},b_{y,t}^{\rho_{t-1}},b_{\mu,t}^{\rho_{t-1}}\}$.
  By equating the coefficients of $\{x_{t-1},y_{t-1},z_{t-1}\}$ with
  $\{a_{x,t}^{\rho_{t-1}},a_{y,t}^{\rho_{t-1}},a_{\mu,t}^{\rho_{t-1}}\}$,
  respectively, we can obtain \eqref{eq:simple-start2},
  \eqref{eq:recursion-ay2} and
  \ref{eq:recursion-amu2}. \eqref{eq:recursion-by-bmu2} can be derived by
  considering \eqref{eq:singlestageproblemV2} in the same
  way. \eqref{eq:secondordercondition2} corresponds to the second order
  conditions for the two maximization problems.
\end{proof}

\begin{th:linear-pi-psi}
  If $U_t$ is TQD, $\psi_t$ is linear, and $\hat{\pi}_t$ is
  linear, then there exists a linear $\pi_t$ such that
\[
\pi_t(x_{t-1}, y_{t-1},\phi_{t-1}) 
\in \argmax_{u_t}  
\left(F_{u_t}^{(\psi_t,\hat{\pi}_t)} U_t\right)(x_{t-1}, y_{t-1},\phi_{t-1}),
\]
for all $x_{t-1}$, $y_{t-1}$, and Gaussian $\phi_{t-1}$, so long as
the optimization problem is bounded.  Similarly, if $V_t$ is AQD and
$\pi_t$ is linear then there exists a linear $\psi_t$ such that
\[
\psi_t(y_{t-1},\phi_{t-1}) 
\in \argmax_{v_t} 
\left(G_{v_t}^{\pi_t} V_t\right)(y_{t-1},\phi_{t-1}),
\]
for all $y_{t-1}$ and Gaussian $\phi_{t-1}$, so long as the
optimization problem is bounded.
\end{th:linear-pi-psi}
\begin{proof}
Suppose that
\begin{gather*}
\begin{split}
U_t (x_t,y_t,\phi_t) &=
-\lambda\bigg(\tfrac{1}{2} c^{\rho_t}_{xx,t} x_t^2 + \tfrac{1}{2} c^{\rho_t}_{yy,t} y_t^2
 + \tfrac{1}{2} c^{\rho_t}_{\mu\mu,t} \mu_t^2 \\
&\quad\quad\quad
 + c^{\rho_t}_{xy,t} x_t y_t + c^{\rho_t}_{x\mu,t} x_t \mu_t + c^{\rho_t}_{y\mu,t} y_t \mu_t 
- \frac{\sigma_\epsilon^2}{\lambda^2} c^{\rho_t}_{0,t}\bigg),\\
\end{split} \\
\hat{\pi}_t(x_{t-1},y_{t-1},\phi_{t-1}) 
= \hat{a}^{\rho_{t-1}}_{x,t} x_{t-1} 
+ \hat{a}^{\rho_{t-1}}_{y,t} y_{t-1} + \hat{a}^{\rho_{t-1}}_{\mu,t} \mu_{t-1}, 
\\
\psi_t(y_{t-1},\phi_{t-1}) = b^{\rho_{t-1}}_{y,t} y_{t-1} 
+ b^{\rho_{t-1}}_{\mu,t} \mu_{t-1}.
\end{gather*}
If the trader takes the action $u_t$, while the arbitrageur uses the
policy $\psi^*_t$ and assumes that the trader uses the policy
$\hat{\pi}_t$, we have
\[
\begin{split}
v_t & =
b^{\rho_{t-1}}_{y,t} y_{t-1} + b^{\rho_{t-1}}_{\mu,t} \mu_{t-1}, \\
x_t & = x_{t-1} + u_t, \\
y_t & = y_{t-1} + b^{\rho_{t-1}}_{y,t} y_{t-1} + b^{\rho_{t-1}}_{\mu,t} \mu_{t-1}.
\end{split}
\]
Using these facts, Theorem~\ref{th:Gaussian}, and
\eqref{eq:expmu-trader1}--\eqref{eq:expmu-trader2} from Lemma~\ref{le:expmu},
we can explicitly compute
\[
\left(F^{(\psi_t, \hat{\pi}_t)}_{u_t} U_t\right)(x_{t-1},y_{t-1},\phi_{t-1})
= \E^{(\psi_t, \hat{\pi}_t)}_{u_t}
\left[\left.\lambda (u_t + v_t) x_{t-1} 
+ U_t(x_t, y_t, \phi_t)\ \right|\  
x_{t-1}, y_{t-1}, \phi_{t-1}\right].
\]
It is easy to see that $\left(F^{(\psi_t, \hat{\pi}_t)}_{u_t}
U_t\right)(x_{t-1},y_{t-1},\phi_{t-1})$ is quadratic in
$u_t$. Moreover, the coefficient of $u^2_t$ is independent of
$\{x_{t-1}, y_{t-1},\mu_{t-1}\}$ while the coefficient of $u_t$ is
linear in $\{x_{t-1},y_{t-1},\mu_{t-1}\}$. Therefore, the
optimizing $u^*_t$ is a linear function of $\{x_{t-1},y_{t-1},\mu_{t-1}\}$, whose
coefficients can be computed by substitution and rearrangement of the 
resulting terms.

Similarly, suppose that
\begin{gather*}
V_t(y_t,\phi_t) 
= - \lambda\left(\tfrac{1}{2} d^{\rho_t}_{yy,t} y_t^2 
+ \tfrac{1}{2} d^{\rho_t}_{\mu\mu,t} \mu_t^2
+ d^{\rho_t}_{y\mu,t} y_t \mu_t
-  \frac{\sigma_\epsilon^2}{\lambda^2} d^{\rho_t}_{0,t}\right),
\\
\pi_t(x_{t-1},y_{t-1},\phi_{t-1})
= a^{\rho_{t-1}}_{x,t} x_{t-1} + a^{\rho_{t-1}}_{y,t} y_{t-1}
+ a^{\rho_{t-1}}_{\mu,t} \mu_{t-1}.
\end{gather*}
If the arbitrageur takes the action $v_t$ and assumes that the trader
uses the policy $\pi_t$, we have
\[
\begin{split}
u_t &
= a^{\rho_{t-1}}_{x,t} x_{t-1} + a^{\rho_{t-1}}_{y,t} y_{t-1}
+ a^{\rho_{t-1}}_{\mu,t} \mu_{t-1},
\\
y_t & = y_{t-1} + v_t.
\end{split}
\]
Using these facts, Theorem~\ref{th:Gaussian}, and
\eqref{eq:expmu-arb1}--\eqref{eq:expmu-arb2} from Lemma~\ref{le:expmu}, we can
explicitly compute
\[
\left(G_{v_t}^{\pi_t} V_t\right)(y_{t-1},\phi_{t-1}) 
= \E^{\pi_t}_{v_t}
\left[\left.\lambda (\pi_t + v_t) y_{t-1} + V_t(y_t,\phi_t) \ \right|\ y_{t-1}, \phi_{t-1}\right].
\]
It is easily checked that $\left(G_{v_t}^{\pi_t}
V_t\right)(y_{t-1},\phi_{t-1}) $ is quadratic in $v_t$. Moreover, the
coefficient of $v^2_t$ is independent of $\{y_{t-1},\mu_{t-1}\}$ while
the coefficient of $v_t$ is linear in
$\{y_{t-1},\mu_{t-1}\}$. Therefore, the optimizing $v^*_t$ is a linear of
$\{y_{t-1},\mu_{t-1}\}$, whose coefficients can be computed by substitution
and rearrangement of the resulting terms.
\end{proof}

\end{document}